\documentclass[aop,seceqn,citesort,MSNbibl,dvips]{arximspdf}
\usepackage{mathrsfs}

\doi{10.1214/10-AOP536}
\volume{38}
\issue{6}
\pubyear{2010}
\firstpage{2103}
\lastpage{2135}

\makeatletter

\newcommand{\varprojlim}{\operatorname{\mathop{\lim}\limits_{\longleftarrow}}}

\newcommand{\A}{\mathbb A}
\newcommand{\F}{\mathbb F}
\newcommand{\Z}{\mathbb Z}
\newcommand{\N}{\mathbb N}
\newcommand{\R}{\mathbb R}
\newcommand{\Q}{\mathcal Q}

\newcommand{\G}{\mathscr G}

\newcommand{\XX}{\mathbf X}
\newcommand{\be}{\beta}
\newcommand{\Ga}{\Gamma}
\newcommand{\la}{\lambda}
\newcommand{\de}{\delta}
\newcommand{\ka}{\varkappa}
\newcommand{\si}{\sigma}
\newcommand{\Om}{\Omega}
\newcommand{\tth}{\theta}
\newcommand{\inv}{\operatorname{inv}}
\newcommand{\Ex}{\operatorname{Ex}}
\newcommand{\supp}{\operatorname{supp}}
\newcommand{\Path}{\operatorname{Path}}
\newcommand{\wt}{\widetilde}
\newcommand{\Sym}{\mathfrak S}
\newcommand{\M}{\mathfrak M}
\newcommand\Pv{P^{(v)}}

\newproclaim{definition}{Definition}[section]
\newtheorem{proposition}[definition]{Proposition}
\newtheorem{lemma}[definition]{Lemma}
\newproclaim{remark}[definition]{Remark}
\newtheorem{corollary}[definition]{Corollary}
\newtheorem{theorem}[definition]{Theorem}

\makeatother

\begin{document}
\begin{frontmatter}

\title{$\lowercase{q}$-exchangeability via quasi-invariance}
\runtitle{$q$-exchangeability via quasi-invariance}

\begin{aug}
\author[A]{\fnms{Alexander} \snm{Gnedin}\corref{}\ead[label=e1]{A.V.Gnedin@math.uu.nl}} and
\author[B]{\fnms{Grigori} \snm{Olshanski}\thanksref{t2}\ead[label=e2]{olsh2007@gmail.com}}
\runauthor{A. Gnedin and G. Olshanski}
\affiliation{Utrecht University and Institute for Information
Transmission Problems and Independent University of Moscow}
\address[A]{Mathematical Institute\\
Utrecht University\\
Postbus 80010\\
3508 TA Utrecht\\
The Netherlands\\
\printead{e1}} %adresu isvedimo komanda gale!
\address[B]{Institute for Information\\
\quad Transmission Problems\\
Bolshoy Karetny 19\\
Moscow 127994\\
and\\
Independent University of Moscow\\
Russia\\
\printead{e2}}
\end{aug}

\thankstext{t2}{Supported by Utrecht University,
RFBR Grant 08-01-00110 and the Project SFB 701 of Bielefeld University.}

% HISTORY:
\received{\smonth{8} \syear{2009}}
\revised{\smonth{1} \syear{2010}}

% ABSTRACT
%
\begin{abstract}
For positive $q\neq1$, the $q$-exchangeability of an infinite random
word is introduced as quasi-invariance under permutations of letters,
with a special cocycle which accounts for inversions in the word. This
framework allows us to extend the $q$-analog of de Finetti's theorem
for binary sequences---see Gnedin and Olshanski [\textit{Electron. J. Combin.} \textbf{16} (2009) R78]---to general
real-valued sequences. In contrast to the classical case of
exchangeability ($q=1$), the order on $\mathbb R$ plays a significant
role for the $q$-analogs. An explicit construction of ergodic
$q$-exchangeable measures involves random shuffling of
$\N=\{1,2,\ldots\}$ by iteration of the geometric choice. Connections
are established with transient Markov chains on $q$-Pascal pyramids and
invariant random flags over the Galois fields.
\end{abstract}

% KEYWORDS
%

\end{frontmatter}

%s1 ###
\section{Introduction}\label{1}

A random word $w=w_1w_2\cdots$ with letters $w_i\in\A$ over some alphabet
$\A\subseteq{\mathbb R}$ is \textit{exchangeable} if swapping the places
of two
neighboring letters $w_i$ and $w_{i+1}$ does not change the
probability. We
shall study the following deformation of this fundamental random symmetry
property. For positive parameter $q\neq1$, we define $w$ to be
\textit{$q$-exchangeable} if, by swapping the places of two neighboring letters $w_i$
and $w_{i+1}$, the probability is multiplied by factor
$q^{\mathrm{sgn}(w_{i+1}-w_i)}$. The intuitive effect of the deformation is that the
arrangement of letters in the word is not completely random, as for exchangeable
sequences, but rather there is a tendency for some monotonic pattern.
To be
definite, we shall focus temporarily on the instance $0<q<1$, in which case,
words with smaller numbers of \textit{inversions} are more likely, the latter
defined as pairs of positions $i<j$ with $w_i>w_j$.

The same definitions apply to finite words. It is well known and easy
to see
that the most general exchangeable word of fixed length $n$ can be
produced by
first choosing an inversion-free word $v_1\leq\cdots\leq v_n$ from an
arbitrary
probability distribution on the space of weakly increasing sequences
over $\A$,
then shuffling the letters by an independent uniformly random
permutation of
${\mathbb N}_n:=\{1,\ldots,n\}$. We will show that the general finitely
$q$-exchangeable word is produced similarly, with the amendment that
the random
permutation should follow the Mallows distribution \cite{Mallows},
which assigns
to each particular permutation $\sigma\dvtx{\mathbb N}_n\to{\mathbb N}_n$
probability proportional to $q^{\inv(\sigma)}$, with $\inv(\sigma)$
being the
number of inversions in $\sigma$.

The analogy between the above exchangeable and $q$-exchangeable representations
does not extend to infinite words. According to de Finetti's theorem
\cite{Al,K}, an infinite exchangeable $w$ satisfies the strong law of
large numbers: for
fixed $B\subset\A$, the proportion of letters $w_i\in B$ in the initial
subword of length $n$ is asymptotic to $\nu(B)$, where $\nu$ is a random
probability measure on $\A$. Conditionally on $\nu$, the random word is
distributed like an i.i.d. sample from $\nu$, so the ergodic
distributions for $w$
are parametrized by probability measures on $\A$. In contrast to that, there
is no principal difference between the general representations of
finite and
infinite $q$-exchangeable words. According to our main result (Theorem
\ref{main}), every infinite $q$-exchangeable $w$ can be produced by choosing
a random sequence $v_1\leq v_2\leq\cdots$ from some arbitrary
distribution on
the space of infinite increasing sequences over $\A$, then shuffling the
letters in the order determined by an independent permutation $\sigma
\dvtx{\mathbb
N}\to{\mathbb N}$ whose distribution is a properly generalized Mallows
distribution on the group $\mathfrak S$ of all bijections of the set
$\N
$. Thus,
every ergodic $q$-exchangeable distribution for $w$ is supported by a single
orbit of the group $\mathfrak S$ acting on $\A^\infty$ by
permutations of
coordinates.

Proving the stated representation of $q$-exchangeable words and
analysis of
the Mallows distribution on $\mathfrak S$ constitute the main contents
of this
paper. For $\A={\mathbb N}_d$, we show that every $q$-exchangeable $w$
can be
encoded into an increasing random walk on the $d$-dimensional lattice with
weighted edges (the $q$-Pascal pyramid); the ergodic measures are
derived, in
this case, by solving a boundary problem via path counting and
asymptotics of
the Gaussian multinomial coefficients.

In our recent paper \cite{GO2}, we observed that every homogeneous random
subspace of an infinite-dimensional space $V$ over a Galois field corresponds
to a $q$-exchangeable sequence over $\A=\{0,1\}$, for $q$ reciprocal
to the
cardinality of the field. Here, homogeneity means \textit{invariance} of the
measure under the natural action on the subspaces of $V$ by the
countable group
of matrices $\mathit{GL}(\infty)=\bigcup_{n\in{\mathbb N}}\mathit{GL}(n)$. In what follows,
we shall
extend this line of research by connecting nonstrict random flags
(sequences of embedded
subspaces in $V$) with $q$-exchangeable words over $\A={\mathbb N}$.

%s2 ###
\section{$q$-exchangeability}\label{2}

In terms of measure theory, $q$-ex\-change\-abi\-lity means
\textit{quasi}-invariance of a probability measure on $\A^\infty$ with
respect to
permutations of an arbitrary finite collection of coordinates, with a special
Radon--Nikodym derivative depending on the altered number of inversions in
the word. To develop this viewpoint, we first recall a general
framework and
some necessary facts from ergodic theory \cite{GS}.

Let $W$ be
a standard Borel space and $G$ be a countable group acting on $W$ on the
left by Borel isomorphisms $T_g\dvtx W\to W$, $g\in G$. Then, $G$ also acts
on the
space of all Borel probability measures on $W$: namely, $T_g$
transforms such
a measure $P$ to $T_gP:=P\circ T_g^{-1}$. We prefer to write this
relation as
$T_g^{-1}P=P\circ T_g$, which means that $(T_g^{-1}P)(X)=P(T_g(X))$ for every
Borel set $X\subseteq W$.

A probability measure $P$ on $W$ is said to be \textit{quasi-invariant} if
$T_g^{-1}P$ is equivalent to $P$ for all $g\in G$, that is, $T_g^{-1}P$
and $P$
have the same null sets.\vspace*{2pt} There then exists a function $\rho(g,w)$ on
$G\times
W$ such that $w\mapsto\rho(g,w)$ is Borel and $T_g^{-1}P=\rho(g,
\cdot
)P$ for
each $g\in G$, that is, $\rho(g, \cdot)$ is the Radon--Nikodym derivative
$dT_g^{-1}P/dP$. The function $\rho$ is unique modulo $P$-null sets and
satisfies
the relation
\[
\rho(gh,w)=\rho(g,T_hw)\rho(h,w),\qquad g,h\in G, w\in W
\]
(again modulo null sets). A function $\rho$ with this property is
called a \textit{multiplicative cocycle}.

Conversely, given a multiplicative cocycle $\rho$, let $\M(\rho)$
denote the
set of all quasi-invariant probability measures on $W$ satisfying the relation
$dT_g^{-1}P/dP=\rho(g, \cdot)$, $g\in G$. The set $\M(\rho)$ has
itself the
structure of a standard Borel space and if $\M(\rho)$ is nonempty, then
it is
convex and has a nonempty subset $\Ex\M(\rho)$ of extreme points. The
set of
extremes $\Ex\M(\rho)$ is also Borel. Moreover, every measure $M\in
\M
(\rho)$
is uniquely representable as a mixture of the extreme measures, meaning that
there exists a unique probability measure $\ka$ on $\Ex\M(\rho)$
such that
\[
M(X)=\int_{\Ex\M(\rho)}P(X)\ka(dP)
\]
for every Borel subset $X\subseteq W$.

Since the generic element of $\M(\rho)$ is a unique mixture of
extremes, it is
important to describe as explicitly as possible the set of extremes
$\Ex\M(\rho)$. A useful criterion is that the extreme measures can be
characterized as ergodic measures from $\M(\rho)$. Recall that a
$G$-quasi-invariant probability measure $P$ on $W$ is \textit{ergodic}
if every
$G$-invariant Borel subset of $W$ has $P$-measure 0 or 1. Since the
group $G$
is countable, the ergodicity is equivalent to the formally stronger condition
that every invariant mod 0 subset has measure 0 or 1.

After these general preliminaries, we focus on a concrete instance. We shall
consider the action of the group $G=\Sym_\infty$ on the infinite
product space
$W=\A^\infty$, where $\Sym_\infty$ is the group of bijections
$\sigma
\dvtx{\mathbb
N}\to{\mathbb N}$ moving only finitely many integers and $\A$ is a
Borel subset
of the ordered space $({\mathbb R},<)$. Although we assume $\A
\subseteq
{\mathbb
R,}$ many considerations of the present paper remain valid for an arbitrary
standard Borel space endowed with a Borel-measurable linear order
(e.g., $\R^k$ with the lexicographic order).

Given a finite word $w=w_1 w_2 \cdots w_n\in\A^n$, let
\[
\inv(w_1\cdots w_n):=\#\{(i,j)\mid1\le i<j\le n, w_i>w_j\}
\]
denote the number of inversions in $w$. For an infinite word $w=w_1 w_2
\cdots
\in\A^\infty$, let
\[
\inv_n(w)=\inv(w_1\cdots w_n)
\]
be the number of inversions in the $n$-truncated word $w_1 \cdots w_n$.

For $w\in\A^\infty$ and $\sigma\in\Sym_\infty$, the difference
$\inv_n(T_\sigma w)-\inv_n(w)$ stabilizes as $n$ becomes so large that
$\si(i)=i$ for all $i\ge n$. We set
%
%e2.1 ###
%
\begin{equation}\label{stable}
c(\si,w)=\mbox{stable value of the difference } \inv_n(T_\sigma
w)-\inv_n(w).
\end{equation}
For instance,
if $\si$ is the elementary transposition of $i$ and $i+1$, then
$T_\sigma w$
differs from $w$ only by transposition of the adjacent letters $w_i$
and $w_{i+1}$, and then $c(\si,w)$ equals $1$, $-1$ or $0$, depending
on whether
$w_i<w_{i+1}$, $w_i>w_{i+1}$ or $w_i=w_{i+1}$, respectively.

The function $c(\si,w)$ is an \textit{additive cocycle} in the sense that
\[
c(\si\tau,w)=c(\si, T_\tau w)+c(\tau,w),\qquad \si, \tau\in
\Sym_\infty.
\]
Equivalently, for $q>0$,
%
%e2.2 ###
%
\begin{equation}\label{rho-q}
\rho_q(\si,w):=q^{c(\si,w)}
\end{equation}
is a multiplicative cocycle.
In accordance with the terminology of ergodic theory,
the additive cocycle $c=\log_q \rho_q$ may be also called
the ``modular function.''

Our considerations are based on the following definition.
\begin{definition}\label{2.A} For fixed $q>0$,
a Borel probability measure $P$ on $\A^\infty$ is called
\textit{$q$-exchangeable} if $P$ is quasi-invariant with respect to the
action of the
group~$\Sym_\infty$, with the multiplicative cocycle given by
(\ref{rho-q}).
\end{definition}

Note that it is enough to require that (\ref{rho-q}) holds for the elementary
transpositions because these permutations generate the group $\Sym
_\infty$.
Thus, Definition \ref{2.A} is equivalent to the definition of
$q$-exchangeability
given in the \hyperref[1]{Introduction}.
In the special case $q=1$, the order on $\A$ plays no role, as the cocycle
$\rho_q$ is identically equal to 1, and so our definition becomes
conventional exchangeability.

It is important to understand how $q$-exchangeability behaves under
transformations. For $f\dvtx\A\to{\mathbb B,}$ let $f^\infty$ denote the induced
mapping $\A^\infty\to{\mathbb B}^\infty$ which replaces each letter
$w_i$ in a
word by $f(w_i)$. First, consider the identity mapping from $(\A,<)$ to
$(\A,>)$.
\begin{proposition}\label{reverse}
If $P$ is a $q$-exchangeable measure on the space of words over
$(\A,<)$, then $P$
is $q^{-1}$-exchangeable with respect to $(\A,>)$, that is,
when the order on the basic space is reversed.
\end{proposition}
\begin{pf}
The claim is easily checked for the elementary transpositions which
swap $i$ and $i+1$.
\end{pf}

It is obvious that if $f$ is an injective morphism of ordered Borel spaces,
then $f^\infty$ sends one $q$-exchangeable measure to another $q$-exchangeable
measure. This applies, in particular, to $\A\subseteq\R$ and a strictly
increasing function $f\dvtx\A\to\R$. It is less obvious that
$q$-exchangeability is
preserved by arbitrary monotone transformations.
\begin{proposition}\label{monotone}
Let $\A$ and ${\mathbb B}$ be Borel subsets of $\R$. Suppose $f\dvtx\A
\to
{\mathbb
B}$ is weakly increasing, that is, $a<b$ implies $f(a)\le f(b)$. The
induced Borel map $f^\infty\dvtx\A^\infty\to{\mathbb B}^\infty$ then
preserves $q$-exchangeability.
\end{proposition}

This proposition will be reduced to its restricted version involving finite
random words and a finite alphabet $\A$ (see Proposition
\ref{f-monotone} below). In the case $q=1$, the assertion becomes a
familiar property of
exchangeability, one which holds for arbitrary Borel $f$.

Definition \ref{2.A} has a straightforward counterpart for \textit{finite}
random words $w\in\A^n$. Let $\Sym_n$ denote the group of
permutations of
$\N_n$. We say that a probability measure $P_n$ on $\A^n$ is \textit{finitely
$q$-ex\-change\-able} if, for each $\si\in\Sym_n$, the measure
$T_\si
^{-1} P_n$
is equivalent to $P_n$ and the Radon--Nikodym derivative $d T_\si^{-1}P_n/d
P_n$ is given by the function $q^{\inv(T_\sigma w)-\inv(w)}$. If $\A$
is finite
or countable, then $P_n$ is purely atomic and this condition means
that, for
$w=w_1\cdots w_n\in\A^n$,
%
%e2.3 ###
%
\begin{equation}\label{fin-eq}
P_n(T_\si w)=q^{\inv(T_\sigma w)-\inv(w)}P_n(w),\qquad \si\in\Sym
_n .
\end{equation}

Consider the canonical projection $\A^\infty\to\A^n$ assigning to
an infinite
word $w=w_1 w_2\cdots$ its $n$-truncation $w_1\cdots w_n$,
$n=1,2,\ldots.$ Given
a probability measure $P$ on $\A^\infty$, let $P_n$ stand for the push-forward
of $P$ under the projection. The following result follows easily from
the definitions.
\begin{lemma}\label{fin-n}
A probability measure $P$ on $\A^\infty$ is $q$-exchangeable if and
only if
$P_n$ is finitely $q$-exchangeable for every $n=1,2,\ldots.$
\end{lemma}

In principle, the structure of the set of finitely $q$-exchangeable
measures on
$\A^n$ is clear: by finiteness of the group $\Sym_n$, every such
measure is a
unique mixture of the extreme measures and every extreme (i.e.,
ergodic) measure is
supported by a single $\Sym_n$-orbit in $\A^n$. Moreover, every $\Sym
_n$-orbit
carries a unique $q$-exchangeable probability measure, hence the extreme
measures are in bijective correspondence with the set of $\Sym
_n$-orbits in
$\A^n$. Each $\Sym_n$-orbit in $\A^n$ contains exactly one word
$v_1\cdots v_n\in\A^n$ which is \textit{inversion-free}, that is, which
satisfies $v_1\le\cdots\le
v_n$. Thus, the collection of inversion-free words of length $n$ parametrizes
the orbits of $\Sym_n$ and all extreme finitely $q$-exchangeable
measures on
$\A^n$.

We can now state a simplified version of Proposition \ref{monotone}.
\begin{proposition}\label{f-monotone}
Let $\A$ and $\mathbb B$ be finite ordered alphabets and let $f\dvtx \A\to
{\mathbb
B}$ be a weakly increasing map. The induced map $f^n\dvtx\A^n\to{\mathbb B}^n$
then preserves the finite $q$-exchangeability of measures.
\end{proposition}

We first show how to deduce Proposition \ref{monotone} from Proposition
\ref{f-monotone}. To this end, let $\A$, ${\mathbb B}$ and $f$ be as required
in Proposition \ref{monotone}. Furthermore, let $P$ be a $q$-exchangeable
probability measure on $\A^\infty$ and $f^\infty(P)$ be its
push-forward under~$f^\infty$. Observe that $(f^\infty(P))_n=f^n(P_n)$ for all
$n=1,2,\ldots.$ By
virtue of Lemma \ref{fin-n}, it suffices to prove that if a measure $P_n$
on $\A^n$ is finitely $q$-exchangeable, then so is its push-forward $f^n(P_n)$.
This, in turn, shows that it suffices to inspect the particular case of extreme
$P_n$. As pointed out above, every extreme measure $P_n$ is
concentrated on a
single $\Sym_n$-orbit so that $P_n$ actually lives on words from a finite
alphabet. This provides the desired reduction to Proposition~\ref{f-monotone}.
\begin{pf*}{Proof of Proposition \ref{f-monotone}}
Let $P_n$ be a finitely $q$-exchangeable measure on $\A^n$ and $\wt
P_n=f^n(P_n)$ its push-forward on ${\mathbb B}^n$. Since the alphabets
are finite, the measures are
purely atomic, supported by finite sets, so we may deal with probabilities
of individual words.

It suffices to prove that for every word $u\in{\mathbb B}^n$ and every
elementary
transposition $\si=(i,i+1)$, one has
\[
\wt P_n(u^*)=q^{\inv(u^*)-\inv(u)}\wt P_n(u),\qquad u^*:=T_\si u.
\]
Let us fix $u$ and $i$. There are three possible cases: $u_i=u_{i+1}$,
$u_i<u_{i+1}$ and $u_i>u_{i+1}$. In the first case, $u^*=u$ and the desired
relation is trivial. By symmetry between the second and third cases, it
suffices to examine one of them, say, the second case. Then,
$\inv(u^*)-\inv(u)=1$. Consider the inverse images $X=(f^n)^{-1}(u)$ and
$X^*=(f^n)^{-1}(u^*)$. We then have $\wt P_n(u)=P_n(X)$ and $\wt
P_n(u^*)=P_n(X^*)$. Thus, we are reduced to showing that
\[
P_n(X^*)=q P_n(X).
\]

Since $f$ is weakly increasing, $u_i<u_{i+1}$ implies that
$w_i<w_{i+1}$ for
every \mbox{$w\in X$}, hence $P(T_\si w)=qP(w)$. It remains to note that
the transformation $T_\si\dvtx\A^n\to\A^n$ maps $X$ bijectively onto $X^*$.
This concludes the proof.
\end{pf*}

Another proof will be given at the end of Section \ref{3}.
\begin{proposition}\label{extreme}
Let $f\dvtx\A\to{\mathbb B}$ be as in Proposition \ref{monotone}. If a
probability measure $P$ on $\A^\infty$ is $q$-exchangeable and extreme,
then so is its push-forward $f^\infty(P)$.
\end{proposition}
\begin{pf}
By Proposition \ref{monotone}, $f^\infty(P)$ is $q$-exchangeable,
hence quasi-in\-variant under the action of $\Sym_\infty$. Obviously,
the map
$f^\infty$ commutes with that action. Recall that extremality of
quasi-invariant measures is equivalent to their ergodicity,
so it suffices to show that $f^\infty(P)$ is ergodic if $P$ is such,
but this follows straightforwardly from the definitions.
\end{pf}

%s3 ###
\section{The finite $q$-shuffle}\label{3}

We fix a positive parameter $q$ (later, we will assume that $0<q<1$).
For a finite\vadjust{\goodbreak}
permutation $\si\in\Sym_n$, we denote by $\inv(\si)$ the number of
inversions,
meaning the number of inversions in the permutation word
$\si(1)\cdots\si(n)$. It is well known that
\[
\sum_{\si\in\Sym_n} q^{\inv(\si)}=[n]_q! ,
\]
where
\[
[n]_q!:=[1]_q [2]_q\cdots[n]_q,\qquad [n]_q:=\sum_{i=0}^{n-1}q^i
\]
[this is a particular case of formula (\ref{MacMahon}) below].
\begin{definition}\label{shuffle1}
For $n=1,2,\ldots,$ the \textit{Mallows measure} $\Q_n$ is the probability
measure on $\Sym_n$ defined by
\[
\Q_n(\si)=\frac{q^{\inv(\si)}}{[n]_q!} ,\qquad \si\in\Sym_n.
\]
\end{definition}

The Mallows measure and its relatives, introduced in \cite{Mallows},
have been studied in statistics in the context of ranking problems; see
\cite{DiaconisRam,Benjamini} for connections with card shuffling and exclusion
processes, and \cite{Starr} for a scaling limit of $\Q_n$.

If $q=1$, then $\Q_n$ is just the uniform measure on $\Sym_n$. Thus,
for general $q>0$,
$\Q_n$ may be viewed as a deformation of the uniform measure.

The Mallows measure is the unique finitely $q$-exchangeable measure supported
by the set of permutation words of length $n$, that is, corresponding
to the
inversion-free word $1\,2 \cdots n$.

The measure $\Q_n$ can be characterized by means of an important independence
property partially mentioned in \cite{Mallows} (at the top of
\cite{Mallows}, page 125, substitute $q^{-1/2}$ for Mallows' $\phi$).
First, we
need more notation. For $n=1,2,\ldots,$ we denote by $G_{q,n}$ the $n$-truncated
geometric distribution on $\N_n=\{1,\ldots,n\}$ with parameter $q$:
\[
G_{q,n}(i)=\frac{q^{i-1}}{[n]_q} ,\qquad i\in\N_n.
\]
For permutation $\si\in\Sym_n$, written as the word
$\si(1)\cdots\si(n)$, define \textit{backward ranks}
%
%e3.1 ###
%
\begin{equation}\label{backrank}
\be_j=\be_j(\si):=\#\{i\le j \mid\sigma(i)\le\sigma(j)\},\qquad
j=1,\ldots,n.
\end{equation}
For instance, the permutation word $1324$ has $\be_1=1$, $\be_2=2$,
$\be_3=2$,
$\be_4=4$.
The correspondence $\sigma\mapsto(\beta_1(\si),\ldots,\beta_n(\si
))$ is a
well-known bijection between $\Sym_n$ and the Cartesian product
$\N_1\times\cdots\times\N_n$.
\begin{proposition}\label{ranks-fin}
The Mallows measure $\Q_n$ is the unique measure on $\Sym_n$ under
which the
backward ranks are independent, with each variable
$j-\beta_j+1$ distributed according to $G_{q,j}$.
\end{proposition}
\begin{pf}
Decompose the number of inversions as $\inv(\sigma)=\sum_{j=1}^n
(j-\beta_j)$
and multiply probabilities of the truncated geometric distribution to
see that
$\Q_n$ coincides with the product measure.
\end{pf}

The following shuffling algorithm is central to our construction of finitely
$q$-exchangeable measures.
The procedure is a variation of ``absorption sampling'' which was
studied under
various guises in
\cite{B,Ke1,Rawlings}.
\begin{definition}\label{shuffle2}
Given an arbitrary finite word $v_1\cdots v_n$, its \textit{$q$-shuffle}
is the
random word $w_1\cdots w_n$ obtained by a random permutation of the letters
$v_1,\ldots,v_n$, determined by the following $n$-step algorithm (not
to be
confused with the notion of $a$-shuffle with integer parameter $a$; see
\cite{BD,St,GO1}).

Let $\xi_1,\ldots,\xi_n$ be independent random variables with $\xi
_j$ having
distribution $G_{q,n-j+1}$.

At step 1, take for $w_1$ the $\xi_1$th letter from the word
$v^{(1)}:=v_1\cdots
v_n$. Then, remove the letter $v_{\xi_1}$ from $v^{(1)}$ and denote by
$v^{(2)}$
the resulting word of length $n-1$. Iterate. So, at each following step
$m=2,\ldots,n,$ there is a word $v^{(m)}$ which was derived from the
initial word
by deleting some $m-1$ letters, a new letter $w_m=v^{(m)}_{\xi_m}$ is then
chosen and, if $m<n$, the word $v^{(m+1)}$ is obtained by removing this letter
from $v^{(m)}$.
\end{definition}
\begin{proposition}\label{shuffle3}
Let $v=v_1\cdots v_n$ be an inversion-free word on the ordered alphabet
$\A$,
so $v_1\le\cdots\le v_n$. Let $w$ be the random word obtained from $v$
by the
$q$-shuffle algorithm and let $P_n$ be the distribution of $w$ which is a
probability measure concentrated on the $\Sym_n$-orbit of $v$. Then,
$P_n$ is
finitely $q$-exchangeable.
\end{proposition}
\begin{pf}
First, observe that the probability $P_n(w)$ of any word $w$ from the
$\Sym_n$-orbit of $v$ is strictly positive. By the very definition of
finite $q$-exchangeability, it suffices to prove that if $\si$ is an elementary
transposition\vspace*{1pt} $(i,i+1)$, $i=1,\ldots,n-1$, then the ratio $P_n(T_\si
(w))/P_n(w)$
equals $q$, $q^{-1}$ or 1, depending on whether $w_i<w_{i+1}$,
$w_i>w_{i+1}$ or
$w_i=w_{i+1}$, respectively.
The latter case being trivial, we may assume, by symmetry, that
$w_i<w_{i+1}$.

For $w_1<w_2$, suppose that a word starts with $w_1w_2$ and examine the
transposition $\sigma=(1,2)$, which swaps $w_1$ and $w_2$. Let $I$ and $J$
denote the sets of indices $i$ and $j$ for which $v_i=w_1$ and $v_j=w_2$,
respectively. If the $q$-shuffle algorithm results in the word $w$,
then the
first chosen letter is $v_i$ for some $i\in I$ and the second chosen
letter is
$v_j$ for some $j\in J$. Likewise, if the resulting word starts with $w_2w_1$,
then we have to choose first $v_j$ with some $j\in J$ and afterward $v_i$
with some $i\in I$. Let $P_{v_iv_j}$ and $P_{v_jv_i}$ stand for the
corresponding probabilities.

If we fix $i\in I$ and $j\in J$, then the word $v^{(3)}$ obtained from the
initial word $v$ at the third step of the algorithm does not depend on the
order in which $v_i$ and $v_j$ were chosen. Thus, it suffices to prove that
$P_{v_iv_j}/P_{v_jv_i}=1/q$.

The probabilities in question are easily computed. Note that $i<j$ because
\mbox{$v_i<v_j$}. It follows that
\[
P_{v_iv_j}=G_{q,n}(i)G_{q,n-1}(j-1)=\frac{q^{i+j-3}}{[n]_q[n-1]_q}
\]
because, after the first step, the letter $v_j$ acquires the number
$j-1$. On the
other hand,
\[
P_{v_jv_i}=G_{q,n}(j)G_{q,n-1}(i)=\frac{q^{i+j-2}}{[n]_q[n-1]_q},
\]
because now the position of the second letter does not change after the first
step. Therefore, the ratio in question is indeed equal to $1/q$.

Finally, transpositions $\si=(i,i+1)$ with $i=2,3,\ldots$ are handled in
the same way,
the key point being that each of the words $v^{(2)}$, $v^{(3)},\ldots$
is inversion-free.
\end{pf}
\begin{remark}
Note that the claim of Proposition \ref{shuffle3} fails if one drops the
assumption that $v$ is inversion-free. For instance, if $v_1\ge\cdots
\ge v_n$,
then the resulting probability measure on the orbit will be
$q^{-1}$-exchangeable and hence not $q$-exchangeable, except the
trivial cases
where $v_1=\cdots=v_n$ or $q=1$.
\end{remark}

The connection between Definitions \ref{shuffle1} and \ref{shuffle2} is
established by the following
result.
\begin{corollary}\label{shuffle-Qn}
The $q$-shuffle, as introduced in Definition \ref{shuffle2},
coincides with the
action of the random permutation $\si\in\Sym_n$, distributed according
to the
Mallows measure $\Q_n$.
\end{corollary}
\begin{pf}
As seen from the description of the $q$-shuffle, it actually acts on
the positions of the letters rather than on the letters themselves.
Thus, it is given
by the action of the random permutation $\si\in\Sym_n$, distributed
according to
some probability measure $\Q'_n$ on $\Sym_n$, which does not depend on
the word
to be $q$-shuffled. Let us identify permutations $\si\in\Sym_n$ with the
corresponding permutation words $\si(1)\cdots\si(n)$. Then, $\Q'_n$
can be
characterized as the outcome of $q$-shuffling the inversion-free word
$v=1\cdot2\cdots n$. By Proposition \ref{shuffle3}, $\Q'_n$ is a finitely
$q$-exchangeable probability measure concentrated on the $\Sym_n$-orbit
of $v$.
Such a measure is unique and the orbit can be identified with the group
$\Sym_n$ itself. On the other hand, $\Q_n$ is $q$-exchangeable, thus
$\Q'_n=\Q_n$.
\end{pf}

As yet another application of Proposition \ref{shuffle3}, we obtain an
alternative proof of Proposition \ref{f-monotone}.
\begin{pf*}{Second proof of Proposition \ref{f-monotone}}
We will show that if $P_n$ is an extreme $q$-exchangeable measure on
$\A
^n$, then
so is $f^n(P_n)$. This will imply the claim of the proposition.

By Proposition \ref{shuffle3}, $P_n$ is obtained by the $q$-shuffle
applied to an
inversion-free word $v\in\A^n$. Therefore, the same holds for the measure
$f^n(P_n)$ and the word $f(v):=f(v_1)\cdots f(v_n)$ because the
$q$-shuffle commutes
with the map $f^n$. Since $f$ is weakly increasing, the word $f(v)$ is
inversion-free. Again applying Proposition~\ref{shuffle3}, we get the desired
result.
\end{pf*}

%s4 ###
\section{The infinite $q$-shuffle and statement of the main
result}\label{4}

The above discussion of finite $q$-exchangeability can be summarized as
follows: the extreme finitely $q$-exchangeable probability measures are
parameterized by finite inversion-free words and can be obtained by application
of the $q$-shuffle procedure to these words. Our aim now is to find a
counterpart
of this result for measures on infinite words.
As in Section \ref{2}, we are dealing with an ordered alphabet $(\A
,<)$, where
$\A$ is a Borel subset of $\R$.
Thus far, the parameter $q$ has been an arbitrary positive number, but:

$\bullet$ throughout the rest of the paper we will assume
$0<q<1$.

By Proposition \ref{reverse}, this restriction does not lead to a loss of
generality because the case $q>1$ is reduced to the case $q<1$ by inverting
the order on the alphabet.

Let $\N=\{1,2,\ldots\}$ and let $G_q$ be the geometric distribution on
$\N$ with
parameter~$q$:
\[
G_q(i)=(1-q)q^{i-1},\qquad i\in\N.
\]

\begin{definition}\label{shuffle4}
Let $v=v_1v_2\cdots\in\A^\infty$ be an arbitrary infinite word. The
\textit{infinite $q$-shuffle} of $v$ is the infinite random word
$w=w_1w_2\cdots$
produced by the algorithm similar to that in Definition \ref{shuffle2}. The
only changes are: (i) the independent variables with varying truncated
geometric distributions should be replaced by the independent variables
$\xi_1,\xi_2,\ldots$ with the same geometric distribution $G_q$; (ii) the
number of steps becomes infinite.
\end{definition}

Although the infinite $q$-shuffle involves countably many steps, the
first $n$
letters in the output word $w$ are specified after $n$ steps of the algorithm.
This shows, in particular, that the law of the random word $w$ is well defined
as a Borel probability measure on $\A^\infty$.
\begin{lemma}\label{permut}
The output random word $w$ is a random permutation of the letters of
the input
word $v$. That is, all letters of $v$ appear in $w$ with probability
$1$.
\end{lemma}
\begin{pf}
The probability that the first letter $v_1$ will not be chosen in the
first $m$
steps of the algorithm is equal to $q^m$. As $m\to\infty$, this
quantity goes
to 0 so that $v_1$ will appear in $w$ with probability 1. Iterating this
argument, we arrive at the same conclusion for all other letters.
\end{pf}

As above, we say that an infinite word $v\in\A^\infty$ is
\textit{inversion-free} if it has no inversions, that is, if $v_1\le v_2\le
\cdots.$
\begin{proposition}\label{shuffle5}
If $v\in\A^\infty$ is an inversion-free word, then its $q$-shuffle
produces a
$q$-exchangeable Borel probability measure on $\A^\infty$.
\end{proposition}
\begin{pf} Let $P^{(v)}$ denote the measure in question. For any
$n=1,2,\ldots,$
let $P_n^{(v)}$ be the $n$th marginal measure of $P$, as in Lemma \ref{fin-n}.
The same argument as in the proof of Proposition \ref{shuffle3} shows
that each
of the measures $P_n^{(v)}$ is $q$-exchangeable. Consequently, by
virtue of
Lemma \ref{fin-n}, $P^{(v)}$ is also $q$-exchangeable.
%\rightqed
\end{pf}

Let $\Sym$ stand for the set of \textit{all} permutations (i.e.,
bijections) of
the set $\N$. We will often identify permutations $\si\in\Sym$ with the
corresponding infinite words $\si(1)\si(2)\cdots\in\N^\infty$. In this
way, we
get an embedding $\Sym\hookrightarrow\N^\infty$. It is easy to
check that
$\Sym$ is a Borel subset of $\N^\infty$ so that one can speak about Borel
measures on $\Sym$.

On the other hand, $\Sym$ is a group containing $\Sym_\infty$ as a proper
subgroup. The group $\Sym$ acts on $\A^\infty$ in the same way as
$\Sym_\infty$ does. Namely, if $\si\in\Sym$ and $w\in\A^\infty
$, then
$(T_\si
w)_i=w_{\si^{-1}(i)}$.
\begin{definition}\label{shuffle6}
By virtue of Proposition \ref{shuffle5} and Lemma \ref{permut}, an
application of
the infinite $q$-shuffle to the inversion-free word $v=1 \cdot2 \cdots\in
\N
^\infty$
produces a $q$-exchangeable Borel probability measure on $\N^\infty$,
which is
concentrated on the group $\Sym$. We call this measure the \textit{Mallows
measure on} $\Sym$ and denote it $\Q$.
\end{definition}
\begin{remark}\label{inversion}
In accordance with our definition of the action of permutations on
words, the
permutation word $\si(1)\si(2)\cdots$ corresponding to an element
$\si
\in\Sym$
coincides with $T_{\si^{-1}}(1\cdot2\cdots)$ and not with $T_\si(1\cdot2\cdots
)$. It
follows that the infinite $q$-shuffle of any infinite word coincides
with the
action on it by the random permutation $T_{\si}$ with $\si\in\Sym$
distributed
according to the push-forward of $\Q$ under the inversion map $\si
\mapsto
\si^{-1}$. However, as will be shown in the \hyperref[10]{Appendix}, $\Q$ is actually
preserved by this map, so we may simply choose random $\si$,
itself distributed according to the Mallows measure~$\Q$.
\end{remark}

Given a word $v\in\A^\infty$, its \textit{support}, denoted $\supp
(v)$, is
the subset of $\A$ comprised of all distinct letters that appear in
$v$, without
regard to their multiplicities. If no assumption on $v$ is made, then
$\supp(v)$ may
be any finite or countable subset of $\R$ and the letters from $\supp
(v)$ may
enter $v$ with arbitrary multiplicities, finite or infinite. This is
not the
case, however, if $v$ is inversion-free, as demonstrated by the following,
evident, proposition.

\begin{proposition}\label{word-v}
The inversion-free words $v\in\R^\infty$ belong to one of the following
two types,
depending on whether the support $\supp(v)$ is finite or infinite:

\begin{enumerate}[(II)]
\item[(I)] The finite type: $\supp(v)$ is a finite set
$a_1<\cdots<a_d$. Then, for each $i=1,\ldots,d-1$, the letter $a_i$
enters $v$
with a finite nonzero multiplicity $l_{a_i}$, while the last letter
$a_d$ has
infinite multiplicity and
\[
v=\underbrace{a_1\cdots a_1}_{l_{a_1}}\cdots\underbrace
{a_{d-1}\cdots
a_{d-1}}_{l_{a_{d-1}}} \underbrace{a_d a_d\cdots}_{l_{a_d}=\infty}.
\]

\item[(II)] The infinite type: $\supp(v)$ is a countable set
$a_1<a_2<\cdots.$ Then, for each $i=1,2,\ldots,$ the letter $a_i$ enters
$v$ with
a finite nonzero multiplicity $l_{a_i}$ and
\[
v=\underbrace{a_1\cdots a_1}_{l_{a_1}} \underbrace{a_2\cdots
a_2}_{l_{a_2}}\cdots.
\]
\end{enumerate}
For both types, the finite multiplicities $l_{a_i}$ may take arbitrary
positive integer values.
\end{proposition}

For an inversion-free word $v\in\R^\infty$, let $\Om^{(v)}$ denote its
$\Sym$-orbit, $\Om^{(v)}:=\{T_\si v\mid\si\in\Sym\}$, which is a
Borel subset
in $\R^\infty$. By the definition, the measure $P^{(v)}$ is concentrated
on $\Om^{(v)}$.
\begin{remark}\label{orbit}
If $\supp(v)$ is finite, then $\Om^{(v)}$ coincides with the $\Sym
_\infty$-orbit
of $v$ and hence is countable [except when $\supp(v)$ is a singleton].
Therefore, in this case, the measure $P^{(v)}$ is purely atomic: for
$w\in\Om^{(v)}$, $P^{(v)}(w)$ is proportional to $q^{\inv(w)}$. Note
that, here,
$\inv(w)$, the total number of inversions in $w$, is finite. Moreover, the
number
\[
\mathcal I^{(v)}(k):=\# \bigl\{w\in\Om^{(v)}\mid\inv(w)=k\bigr\}
\]
has polynomial growth in $k$ as $k\to\infty$ so that the series $\sum_k
\mathcal I^{(v)}(k)q^k$ converges, which explains why the measure exists.
(Note that in the situation of the conventional de Finetti theorem,
there are
no finite invariant measures supported by a nontrivial $\Sym_\infty
$-orbit.) In
contrast to that, if $\supp(v)$ is infinite, then $\Om^{(v)}$ has the
cardinality of the
continuum and the measure $P^{(v)}$ is diffuse.
\end{remark}

We are now in a position to state the main result of the paper.
\begin{theorem}\label{main}
Let $\A$ be an arbitrary Borel subset of $\R$ with order inherited from
$\R$. The extreme $q$-exchangeable Borel probability measures on $\A
^\infty$
are parametrized by the infinite inversion-free words $v$ with support
contained in $\A$. The measure $P^{(v)}$ corresponding to such a word
$v$ is
obtained by application of the infinite $q$-shuffle to $v$, as
described in
Proposition \ref{shuffle5}.
\end{theorem}

Observe that the orbits $\Om^{(v)}$ with different $v$'s are pairwise disjoint.
It follows that Theorem \ref{main} is reduced to the following,
seemingly weaker,
claim.
\begin{proposition}\label{weak}
For $\A$ as in Theorem \ref{main}, the extreme $q$-exchangeable
measures on
$\A^\infty$ belong to the family of measures $\{P^{(v)}\}$, where $v$ ranges
over the set of inversion-free words in $\A^\infty$.
\end{proposition}

Indeed, combining this proposition with the above observation, we see
that none
of the measures in the family $\{P^{(v)}\}$ can be written as nontrivial
mixtures of
other measures, which implies that each $P^{(v)}$ is extreme.
A proof of Proposition \ref{weak} will be given below.
\begin{remark}\label{main1}
Given an element $\tau\in\Sym$, let $\wt\tau\in\N^\infty$
denote the
corresponding permutation word, $\wt\tau=\tau(1)\tau(2)\cdots.$
The Mallows
measure $\Q$ (Definition \ref{shuffle6}) can be characterized as the only
probability measure on the group $\Sym$, which is quasi-invariant
under the
right shifts $\tau\mapsto\tau\si^{-1}$ by elements $\si$ of the subgroup
$\Sym_\infty$, with the cocycle $\rho_q(\si,\wt\tau)$. This
follows from
Theorem \ref{main} and the definition of $\Q$.
\end{remark}

Next, we shall inspect the nature of the random word $w\in\A^\infty$ under
$P^{(v)}$. The sequence of truncations $\varnothing$, $w_1$,
$w_1w_2,\ldots$ has
transition probabilities described in the following proposition. The
notation works as follows: letters $a$, $b$ range over~$\A$; $l_a$ is
the multiplicity of $a$ in $v$, as
above; $u=w_1\cdots w_{n-1}$ is a finite word; $\mu_a(u)$ is the
multiplicity of
$a$ in $u$.
\begin{proposition}\label{trans}
Let $w$ be the infinite random word distributed according to $P^{(v)}$. The
transition probabilities then have the form
%
%e4.1 ###
%
\begin{eqnarray}\label{trans1}
P^{(v)}(u\to ua)&=&q^{\sum_{b<a}(l_b-\mu_b(u))}\bigl(1-q^{l_a-\mu_a(u)}\bigr)
\nonumber\\[-8pt]\\[-8pt]
&=&q^{\sum_{b<a}(l_b-\mu_b(u))}-q^{\sum_{b\le
a}(l_b-\mu_b(u))}.\nonumber
\end{eqnarray}
\end{proposition}
\begin{pf}
First, assume that $n=1$, that is, $u=\varnothing$. The left-hand side of
(\ref{trans1}) is then the probability of $w_1=a$, as in the first step
of the
$q$-shuffling algorithm. The string of $a$'s in $v$ starts from position
$i:=1+\sum_{b<a}l_b$ and ends at position $j:=\sum_{b\le a}l_b$.
Therefore, the
probability in question equals
\[
(1-q)(q^{i-1}+\cdots+q^{j-1})=q^{i-1}(1-q^{j-i+1}).
\]
The same quantity appears in the right-hand side of (\ref{trans1}) when
$u=\varnothing$ because then $\mu_b(u)=0$ for all $b\in\A$.

For $n=2,3,\ldots,$ the argument is exactly the same, taking into
account that
we are dealing with the $n$th step of the algorithm and that the word
$v^{(n)}$ is
inversion-free, with letter multiplicities $l'_b=l_b-\mu_b(u)$.
\end{pf}
\begin{remark}
The following comments are relevant to formula (\ref{trans1}):

1. If $\mu_a(u)=l_a$, then (\ref{trans1}) shows that the transition
$u\to ua$ has
probability zero. This agrees with the fact that if
$l_a<\infty$, then the letter $a$ cannot enter the random word more
than $l_a$
times. In particular, if $l_a=0$ [which means that $a\notin\supp(v)$],
then $a$ never
appears.

2. The transition probability $P^{(v)}(u\to ua)$ depends on $u$ only through
the collection of multiplicities $\{\mu_a(u)\}_{a\in\A}$. That is, it
depends only on the
$\Sym_n$-orbit of~$u$.

3. Recall that the support of $v$ is either of the form $a_1<\cdots
<a_d$ or
$a_1<a_2<\cdots.$ Let us set
\[
x_0(u)=1,\qquad x_i(u)=q^{\sum_{j\le i}(l_{a_j}-\mu_{a_j}(u))},
\]
where $j=1,\ldots,d$ or $j=1,2,\ldots$ for finite or infinite support,
respectively. In this notation, (\ref{trans1}) can be rewritten as
%
%e4.2 ###
%
\begin{equation}\label{trans2}
P^{(v)}(u\to ua_i)=x_{i-1}(u)-x_i(u),\qquad a_i\in\supp(v).
\end{equation}
Now, observe that
\[
1=x_0(u)\ge x_1(u)\ge\cdots\ge x_d(u)=0
\]
or
\[
1=x_0(u)\ge x_1(u)\ge x_2(u)\ge\cdots\ge0 \qquad\mbox{with }
\lim_{i\to\infty}x_i(u)=0
\]
for finite or infinite support, respectively. This makes evident the
fact that
the transition probabilities given by (\ref{trans2}) indeed sum to 1.

4. We have deduced formula (\ref{trans1}) from the $q$-shuffling algorithm.
Conversely, starting from (\ref{trans1}), one can easily recover the
algorithm itself.
\end{remark}

Proposition \ref{trans} describes the measures $P^{(v)}$ via transition
probabilities. The next proposition characterizes $\Pv$ in terms of the
marginal measures $P^{(v)}_n$, which are the joint distributions of the first
$n$ letters. Note that $\Pv_n$ is a purely atomic measure because it is
supported by the words $u=u_1\cdots u_n$ with letters $u_i$ from the
finite or
countable set $\supp(v)$ and the set of all such words is finite or countable.
Thus, we may speak about probabilities $\Pv_n(u)$ of individual words.

We recall some standard $q$-notation. Let
\[
(x;q)_0=1,\qquad (x;q)_k:=\prod_{i=0}^{k-1}(1-xq^i),\qquad k=1,2,\ldots.
\]
Likewise, we define $(x;q^{-1})_k$.
Below, we use the same notation as in Proposition~\ref{trans}.
\begin{proposition}\label{marg}
Let $v\in\R^\infty$ be an inversion-free word and let $u$ be a word of
length $n$ with letters belonging to the support of $v$. We have
%
%e4.3 ###
%
\begin{equation}\label{marg1}
\Pv_n(u)=q^{\inv(u)}q^{-\sum_{b<a}\mu_b(u)\mu_a(u)} \prod_a
(q^{l_a};q^{-1})_{\mu_a(u)}q^{\mu_a(u)\sum_{b<a}l_b},
\end{equation}
where $a$ and $b$ assume values in $\supp(v)$.
\end{proposition}

Note that the product over $a\in\supp(v)$ is actually finite, even if
$\supp(v)$
is infinite. This follows from the fact that $\mu_a(u)=0$ implies that the
corresponding factor equals 1 and that there are only finitely many
$a$'s with
$\mu_a(u)\ne0$.
\begin{pf}
Computing the ratio $\Pv_{n+1}(ua)/\Pv_n(u)$ from (\ref{marg1}), one
sees that
the formula agrees with transition probabilities (\ref{trans1}).
\end{pf}

%s5 ###
\section{The case of a finite alphabet}\label{5}

In this section, we prove Proposition \ref{weak} (and hence Theorem
\ref{main})
for a finite alphabet $\A$ with cardinality $d=\#\A\geq2$. The
simplest case,
$d=2$, was examined in \cite{GO2} and we will apply here the same
method. To be
definite, we take $\A=\N_d$. Following the formalism due to Kerov and
Vershik~\cite{VK}, it is insightful to interpret $q$-exchangeability as a property
of measures on the path space of a graded graph (Bratteli diagram) which
captures the branching of orbits of $\Sym_n$ on $\A^n$ as $n$ varies.

Let $\Z_+=\{0,1,2,\ldots\}$ and consider the $d$-dimensional lattice
$\Z_+^d$.
The lattice points will be denoted by $\la$ or $\mu$. We write lattice
points as
vectors $\lambda=(\lambda_1,\ldots,\lambda_d)$ in the canonical basis
$e_1,\ldots,e_d$ and we call $|\la|=\la_1+\cdots+\la_d$ the
\textit{degree} of
$\lambda$. We write $\mu\prec\la$ if $\mu\ne\la$ and $\la-\mu
\in\Z
_+^d$; in
this case, there is a nondecreasing lattice path connecting $\mu$ with
$\la$.

Each $\lambda$ of degree $n$ corresponds to an inversion-free word,
%
%e5.1 ###
%
\begin{equation}\label{v-lambda}
v(\la)=v_1\cdots v_n
=\underbrace{1\cdots1}_{\la_1} \underbrace{2\cdots2}_{\la_2}
\cdots
\underbrace{d\cdots
d}_{\la_d} ,
\end{equation}
where the letter $a$ does not enter if $\la_a=0$. This correspondence
yields a
bijection between $\Sym_n$-orbits in $\A^n$ and vectors $\lambda\in
\Z
_+^d$ of
degree $n$.
\begin{definition}\label{G}
The \textit{$q$-Pascal pyramid of dimension} $d$, denoted $\Ga(q,d)$,
is the
oriented graph with vertex set $\Z_+^d$ and directed edges
$(\la,\la+e_a)$,
endowed with weights
%
%e5.2 ###
%
\begin{equation}\label{weight}
\operatorname{weight}(\la,\la+e_a):=q^{\la_{a+1}+\cdots+\la_d},\qquad
a\in \N_d.
\end{equation}
\end{definition}

Note that $\operatorname{weight}(\la,\la+e_d)=1$ for any $\la$. The
$n$th level
of the graph consists of the vertices $\la\in\Z_+^d$ with $|\la|=n$.
Level $0$
has a sole \textit{root} vertex $\bar0:=(0,\ldots,0)$. A \textit{standard path}
terminating at $\lambda$ is a lattice path which connects $\bar0$ to
$\lambda$
and is nondecreasing in each coordinate. Similarly, we define an \textit{infinite
standard path} in $\Ga(q,d)$ as an infinite coordinatewise
nondecreasing path
with initial vertex $\bar0$.

Observe that there is a natural bijection between $\A^n$ and standard
paths in
$\Ga(q,d)$ of length $n$. By this bijection, a word $w_1\cdots w_n$ is
mapped to
the path
\begin{eqnarray*}
\mu(\varnothing) &=& \bar0,\qquad \mu(w_1)=e_{w_1}, \qquad
\mu(w_1w_2)=e_{w_1}+e_{w_2}, \ldots, \\
\mu(w_1\cdots w_n)&=&e_{w_1}+\cdots+e_{w_n},
\end{eqnarray*}
where the $a$th coordinate of the terminal vertex is equal to the multiplicity
of the letter $a$ in $w_1\cdots w_n$. For $n=1,2,\ldots,$ the
bijections are
consistent and hence define a bijection between $\A^\infty$ and the
set of
infinite standard paths in $\Ga(q,d)$: under this bijection, $w_n=a$
means that
the $n$th edge of the path connects a vertex $\mu(w_1\cdots w_{n-1})$
of degree
$n-1$ with $\mu+e_a$. Fixing the first $n$ vertices of a standard path
corresponds to a cylinder $[w_1\cdots w_n]\subset\A^\infty$. A measure
$P$ on
$\A^\infty$ translates as a measure on the space of infinite standard paths,
with $P([w_1\cdots w_n])$ being the probability of the corresponding initial
path of length $n$.
\begin{definition}\label{3.H}
The \textit{weight} of a standard path with endpoint $\lambda$ is
defined as the
product of the weights of the edges comprising the path. Let us say
that a
probability measure on the path space of $\Ga(q,d)$ is a \textit{Gibbs
measure}
if, for every~$\lambda$, the conditional measure of a standard path terminating
at $\lambda$ is proportional to the weight of this path (in the
terminology of
Kerov and Vershik \cite{VK}, such a measure is called ``central'').
\end{definition}
\begin{proposition}\label{corr}
For $\A=\N_d$, the $q$-exchangeable measures on $\A^\infty$ correspond
bijectively to the Gibbs measures on the space of infinite standard
paths in
the $q$-Pascal pyramid $\Ga(q,d)$.
\end{proposition}
\begin{pf}
Let $w\in\A^\infty$. Under the correspondence between words and paths,
$q^{\inv_n(w)}$ is equal to the weight of the standard path encoded in
$w_1\cdots w_n$, as seen by induction. Indeed, if the finite word
$w_1\cdots
w_{n-1}$ corresponds to $\lambda$ and $w_n=a$ is appended, then the
number of
inversions increases by
$\inv_{n}(w)-\inv_{n-1}(w)=\lambda_{a+1}+\cdots+\lambda_d$, which
is the same
quantity that appears in (\ref{weight}); we then use the telescoping
representation
\begin{eqnarray*}
\inv_n(w)&=&[\inv_n(w)-\inv_{n-1}(w)]+
[\inv_{n-1}(w)-\inv_{n-2}(w)]+\cdots\\
&&{}+[\inv_1(w)-0].
\end{eqnarray*}

On the other hand, the words in $\A^n$ that correspond to standard
paths with a
given endpoint make up a $\Sym_n$-orbit. Thus, we see that the Gibbs condition
for fixed $n$ is equivalent to finite $q$-exchangeability. Since this
holds for
every $n$, Lem\-ma~\ref{fin-n} allows finite $q$-exchangeability
for $n=1,2,\ldots$ to be translated as the Gibbs property, and conversely.
\end{pf}

We shall now proceed along the lines of \cite{KOO}. Denote by $\Path
(d)$ the
space of all infinite standard paths in $\Ga(q,d)$. With each $\lambda
\in
{\mathbb Z}_+^d$, we associate a unique \textit{elementary} probability measure
supported by the finite set of standard paths with endpoint $\lambda$. This
measure corresponds to an orbital, finitely $q$-exchangeable
probability measure
on $\A^n$. We can understand this measure as a function which assigns to
$\lambda$ value $1$ and assigns to each $\mu\prec\lambda$ the
probability that
a path passes through~$\mu$. The \textit{Martin boundary} of $\Ga(q,d)$ consists
of probability measures on $\Path(d)$ which are representable as weak
limits of these
elementary measures along a sequence of lattice points with
$|\lambda|\to\infty$. We will prove that under the correspondence of
Proposition \ref{corr}, the Martin boundary is exactly the images of the
measures $P^{(v)}$, with $v$ ranging over the set of inversion-free
words in
$\A^\infty$. By the general theory (see \cite{KOO}), the Martin boundary
contains all extreme Gibbs measures, so this will imply Proposition
\ref{weak}.

To determine the boundary, we need to identify all asymptotic regimes for
$\lambda$ which guarantee convergence of the ratios
%
%e5.3 ###
%
\begin{equation}\label{rat}
\frac{\dim(\mu,\lambda)}{\dim(\lambda)},
\end{equation}
where $\dim(\lambda)=\dim(\bar0,\lambda)$ and $\dim(\mu,\lambda
)$ is
equal to
the sum of weights of all nondecreasing lattice paths connecting $\mu$ and
$\lambda$ (the weight of each such path is defined as the product of the
weights of its edges).
We set $\dim(\mu,\la)=0$ if
$\la-\mu\notin{\mathbb Z}_+^d$.
The ratio (\ref{rat}) is the Martin kernel for a certain Markov chain and,
by analogy with the Gibbs formalism
in statistical physics, $\dim\lambda$ may be
called the ``partition function.''

Recall the notation
\[
[0]_q!=1,\qquad [n]_q!=[1]_q[2]_q\cdots[n]_q=\frac
{(q;q)_n}{(1-q)^n},\qquad
n=1,2,\ldots.
\]
For nonnegative integers $n_1,\ldots,n_d$ with $n_1+\cdots+n_d=n$, the number
\[
\left[\matrix{n \cr n_1,\ldots,n_d}\right]_q:=\frac
{[n]_q!}{[n_1]_q!\cdots[n_d]_q!}
=\frac{(q;q)_n}{(q;q)_{n_1}\cdots(q;q)_{n_d}}
\]
is known as the \textit{Gaussian multinomial coefficient}.
\begin{lemma}\label{dim}
We have, for $\lambda=(\lambda_1,\ldots,\lambda_d)$ and $\mu\prec
\lambda$,
\[
\dim(\lambda)=
\left[\matrix{|\la|\cr
\la_1,\ldots,\la_d}\right]_q,\qquad
\dim(\mu,\lambda)=q^{N(\mu,\la)}\dim(\lambda-\mu),
\]
where
\[
N(\mu,\la)= \sum_{b<a}\la_b \mu_a - \sum_{b<a}\mu_b\mu_a.
\]
\end{lemma}
\begin{pf}
Recall that the set of finite standard paths ending at $\la$ is encoded
by the
words $w$ belonging to the $\Sym_{|\la|}$-orbit of the inversion-free word
$v(\la)$, as defined in (\ref{v-lambda}). Let $\{w\}$ stand for the set
of these
words. MacMahon's formula for the generating function for the number of
inversions in permutations of a multiset (see \cite{An}, Theorem 3.6)
says, in
our notation, that
%
%e5.4 ###
%
\begin{equation}\label{MacMahon}
\sum_{\{w\}} q^{\inv(w)}= \left[\matrix{|\la|\cr
\la_1,\ldots,\la_d}
\right]_q.
\end{equation}
This yields the formula for $\dim(\lambda)$.
The formula for $\dim(\mu,\la)$ with
\begin{eqnarray*}
N(\mu,\la)&=&(\la_1-\mu_1)(\mu_2+\cdots+\mu_d)+(\la_2-\mu_2)(\mu
_3+\cdots+\mu_d)+\cdots\\
&&{}+(\la_{d-1}-\mu_{d-1})\mu_d
\end{eqnarray*}
follows
by counting inversions in the corresponding words, which, in turn, is
done by comparing the oriented subgraph
rooted at $\mu$ with the whole graph $\Ga(q,d)$.
\end{pf}

A weakly increasing function $h\dvtx{\mathbb N}_d\to\{0,1,\ldots,\infty\}
$ with
$h(d)=\infty$ will be called a \textit{height function} on $\A={\mathbb
N}_d$. We
also set $h(0):=0$, where appropriate. There is a natural bijection
$h\leftrightarrow v$ between the height functions on $\N_d$ and the
inversion-free words in $\N_d^\infty$,
%
%e5.5 ###
%
\begin{equation}\label{3.f}
v=\underbrace{1\cdots1}_{h(1)} \underbrace{2\cdots2}_{h(2)-h(1)}
\underbrace{3\cdots3}_{h(3)-h(2)} \cdots \underbrace{r\cdots
r}_{h(r)-h(r-1)} \underbrace{r+1\cdot r+1\cdots}_{h(r+1)=\infty} ,
\end{equation}
where, for some $0\le r<d$, each letter $1\le a\leq r$ appears
$h(a)-h(a-1)<\infty$ times (if any), and infinitely many times for $a=r+1$.
\begin{proposition}\label{martin}
The Martin boundary of the graph $\Ga(q,d)$ can be parametrized, in a
natural way, by the height functions on ${\mathbb N}_d$.
\end{proposition}
\begin{pf}
Using the identity
\[
\frac{(q;q)_n}{(q;q)_{n-m}}=(q^n;q^{-1})_m,\qquad n\geq m\geq0,
\]
we derive, from Lemma \ref{dim} for $\mu\prec\lambda$, $m=|\mu|$ and
$n=|\la|$,
that
%
%e5.6 ###
%
\begin{equation}\label{3.a}\quad
\frac{\dim(\mu,\la)}{\dim\la}
=q^{-\sum_{b<a}\mu_b\mu_a} \frac{(q;q)_{n-m}}{(q;q)_n} \prod_{a=1}^d
(q^{\la_a};q^{-1})_{\mu_a}q^{\mu_a\sum_{b\dvtx b<a}\la_b}.
\end{equation}
Observe that the constraint $\mu\prec\la$ can be removed; indeed, if
it is
not satisfied, then $\dim(\mu,\la)=0$ and the right-hand side of
(\ref{3.a})
also vanishes because $(q^{\la_a};q^{-1})_{\mu_a}=0$ for $\la_a<\mu_a$.

Let us rewrite (\ref{3.a}) using the notation
\[
h_\la(a):=\la_1+\cdots+\la_a,\qquad a=1,\ldots,d,\qquad h_\la(0):=0,
\]
in the form
%
%e5.7 ###
%
\begin{eqnarray}\label{3.b}
\frac{\dim(\mu,\la)}{\dim\la}
&=&q^{-\sum_{b<a}\mu_b\mu_a}
\frac{(q;q)_{n-m}}{(q;q)_n}\nonumber\\[-8pt]\\[-8pt]
&&{}\times\prod_{a=1}^d
\bigl(q^{h_\la(a)-h_\la(a-1)};q^{-1}\bigr)_{\mu_a}q^{\mu_ah_\la(a-1)}.\nonumber
\end{eqnarray}
It is now easy to analyze the asymptotics of this expression, assuming that
$\mu$ remains fixed while $\lambda$ varies so that $n=|\la|\to
\infty$.
First, note that
\[
\lim_{n\to\infty}\frac{(q;q)_{n-m}}{(q;q)_n}=\frac{(q;q)_\infty
}{(q;q)_\infty}=1.
\]
Next, observe that
\[
0\le h_\la(1)\le\cdots\le h_\la(d-1)\le h_\la(d)=n.
\]
Passing to a subsequence, we may assume that there exist finite or infinite
limits
\[
\lim_{n\to\infty}h_\la(a)=h(a)\in\Z_+\cup\{+\infty\},\qquad
a=1,\ldots,d .
\]
This means that there exists $0\le r<d$ such that the numbers $h_\la
(1),\ldots,
h_\la(r)$ stabilize for $n$ large enough, $h_\la(a)=h(a)<\infty$ for
$1\leq
a\le r$, while $h_\la(a)\to h(a)=+\infty$ for $a>r$. Note that $h_\la(d)=n$
always goes to infinity so that $h(d)=\infty$ in any case.

Clearly, the product in (\ref{3.b}) up to $a=r$ stabilizes. Next, we have
\[
\bigl(q^{h_\la(r+1)-h_\la(r)};q^{-1}\bigr)_{\mu_{r+1}}q^{\mu_{r+1}h_\la
(r)}\to
q^{\mu_{r+1}h_\la(r)},
\]
because $q^{h_\la(r+1)-h_\la(r)}\to0$. As for the factors with $a>
r+1$, we
have
\[
\bigl(q^{h_\la(a)-h_\la(a-1)};q^{-1}\bigr)_{\mu_a}q^{\mu_ah_\la(a-1)}\to\de
_{\mu_a,0}
\]
with the Kronecker delta in the right-hand side because $h_\la(a-1)\to
\infty$.

We conclude that the convergence $h_\lambda\to h$ implies
%
%e5.8 ###
%
\begin{equation}\label{3.c}
\frac{\dim(\mu,\la)}{\dim\la} \to q^{-\sum_{b<a}\mu_b\mu_a}
\prod_{a=1}^d \bigl(q^{h(a)-h(a-1)};q^{-1}\bigr)_{\mu_a}q^{\mu_ah(a-1)}
\end{equation}
with the convention that $h(0)=0$ and $h(a)-h(a-1)=0$ if
$h(a)=h(a-1)=+\infty$.
Since, for distinct $h$, the limits in (\ref{3.c}) are all distinct,
the Martin boundary can indeed be parameterzed by the height functions.
\end{pf}

Observe that if $h(a)=h(a-1)$, then the limit value (\ref{3.c})
vanishes unless
\mbox{$\mu_a=0$}. Returning to random words $w=w_1w_2\cdots\in\A^\infty$,
this means
that if $h(a)=h(a-1)$, then the letter $a$ does not occur in $w$, with
probability 1.
\begin{proposition}\label{coinc}
Under the correspondence $h\leftrightarrow v$, the measures on $\Path(d)$
afforded by Proposition \ref{martin} correspond exactly to the measures
$\Pv$,
where $v$ ranges over the set of inversion-free words on the alphabet
$\N_d$.
\end{proposition}
\begin{pf}
Fix a height function $h$ and let $\mathcal P$ be the corresponding Gibbs
measure on $\Path(d)$. Next, let $P$ be the measure on $\N_d^\infty$ which
corresponds to $\mathcal P$ via the bijection of Proposition \ref{corr}.
Finally, let $v\in\N_d^\infty$ be the inversion-free word associated
with $h$.
We have to prove that $P=\Pv$. To do this, it suffices to check that
$P_n=\Pv_n$
for all $n$. Let $u\in\N_d^n$. Then, $P_n(u)$ equals $q^{\inv(v)}$
times the
right-hand side of (\ref{3.c}), where we set $\mu_a=\mu_a(u)$.
Comparing with
(\ref{marg1}), we see that this coincides with $\Pv_n(u)$.
\end{pf}

This concludes the proof of Proposition \ref{weak} in the case of a finite
alphabet~$\A$.

%s6 ###
\section{The case $\A=\N$}\label{6}

In this section, we assume that $\A$ is the countable ordered set
$(\N,<)$ of positive integers.
Our aim is to prove, for this case, Proposition \ref{weak}
and hence Theorem \ref{main}.
\begin{definition}\label{height} By a \textit{height function} on $\N$,
we mean a
map $h\dvtx\N\to\Z_+\cup\{+\infty\}$ which is weakly increasing [i.e.,
$h(a)\le
h(b)$ for $a<b$] and satisfies $\lim_{a\to\infty} h(a)=+\infty$. The
set of all
height functions on $\N$ will be denoted $H(\N)$.
\end{definition}

Obviously, setting
\[
l_a=h(a)-h(a-1),\qquad a\in\N,
\]
with the understanding that $h(0)=0$ and $l_a=0$ if
$h(a)=h(a-1)=+\infty
$, we
get a bijection $h\leftrightarrow v$ between $H(\N)$ and the set of all
inversion-free words $v\in\N^\infty$.
\begin{pf*}{Proof of Proposition \ref{weak} for $\A=\N$}
Assume that $P$ is an extreme
$q$-exchangeable measure on $\N^\infty$.
We have to show
that $P=\Pv$ for some $v$. The idea is to reduce this claim to the case
$\A=\N_d$, which was examined in Section \ref{5}, by using Propositions
\ref{monotone} and \ref{extreme}.

For $d=1,2,\ldots$ and $a\in\N$, set $f_d(a)=a\wedge d=\min(a,d)$.
Clearly, this
gives us a weakly increasing map $f_d\dvtx\N\to\N_d$. By Proposition
\ref{extreme},
$f_d^\infty(P)$ is an extreme $q$-exchangeable measure on $\N
_d^\infty
$. By the
results of Section \ref{5}, it coincides with some measure
$P^{(v(d))}$, where
$v(d)\in\N_d^\infty$ is an inversion-free word. Denote by $h_d$ the
corresponding height function on $\N_d$.

Let $w\in\N^\infty$ be the random word with law $P$. For each
$a=1,\ldots,d-1$,
the letter $a$ enters the random word $f_d(w)$ exactly
$h_d(a)-h_d(a-1)$ times,
with probability 1. Since the map $f_d$ does not change the letters
$a=1,\ldots,d-1$, the same holds for the initial random word $w$. This implies
that $h_d(a)=h_{d+1}(a)$ for all $a=1,\ldots, d-1$. Therefore, for every
$a\in\N$, the value $h_d(a)$ stabilizes as $d\to\infty$, starting
from $d=a+1$;
denote by $h(a)$ this stable value. We claim that $h$ is a height
function on
$\N$. Indeed, it is obvious that $h$ weakly increases, so we only have to
check that $h(a)\to\infty$ as $a\to\infty$. If this were not the
case, then
$h(a)$ would assume the same (finite) value for all $a$ large enough.
However, this
would mean that $w$ contained only finitely many letters, each with a
prescribed finite multiplicity $l_a=h(a)-h(a-1)$, which is clearly impossible.
Thus, $h$ should be a height function.

Now, let $v\in\N^\infty$ be the inversion-free word corresponding to
$h$. By
the definition of $h$, we have $f_d^\infty(P)=f_d^\infty(\Pv)$ for all
$d$. Clearly,\vspace*{1pt} this implies $P_n=\Pv_n$ for all $n$, so $P=\Pv$, as
desired.
\end{pf*}
\begin{remark}\label{pascal-infty}
An alternative proof can be based on the notion of the \textit{$q$-Pascal pyramid
of dimension $\infty$}, denoted $\Ga(q,\infty)$, which is the graph
with the
vertex set
\[
\{\la\in\Z_+^\infty\mid\la_1+\la_2+\cdots<+\infty\},
\]
the edges $(\lambda,\lambda+e_a)$, where
\[
e_a=( \underbrace{0,\ldots,0}_{a-1} ,1,0,0,\ldots),\qquad a\in\N,
\]
and the weight $q^{\sum_{b>a}\lambda_b}$ assigned
to the edge
$(\lambda,\lambda+e_a)$. Note
that the sum in the exponent is finite because $|\lambda|:=\sum
_a\lambda
_a$ is
finite, by the definition of $\Ga(q,\infty)$. The $n$th level of
$\Ga(q,\infty)$ consists of vertices with $|\la|=n$.

The graph $\Ga(q,d)$ is embedded in $\Ga(q,\infty)$ as the set of
vertices with
$\la_{b}=0$ for $b>d$. Obviously, $\Ga(q,\infty)=\bigcup_{d\geq1}\Ga
(q,d)$. The
definition of Gibbs measures on the space of standard paths in $\Ga
(q,\infty)$
and the correspondence with $q$-exchangeable measures on $\N^\infty$
straightforwardly extend the definitions from Section \ref{5}. One can
then repeat the
arguments in Proposition \ref{martin} to show that the Martin boundary of
$\Ga(q,\infty)$ consists precisely of the Gibbs measures
corresponding to
measures $\Pv$.
\end{remark}

%s7 ###
\section{The case $\A=\R$}\label{7}

Here, we prove Proposition \ref{weak} and hence Theorem \ref{main} for
$\A=\R$.
This will also cover the seemingly more general case where $\A$ is an arbitrary
Borel subset of $(\R,<)$.

Assume that the measure $P$ on $\R^\infty$
is $q$-exchangeable and extreme. Our aim is to show that there
exists a finite or countable subset $A\subset\R$, of the form
$a_1<\cdots<a_d$ or $a_1<a_2<\cdots,$ such that $P$ is supported by
$A^\infty$. The results of Sections \ref{5} and \ref{6} will then
imply that
$P=\Pv$ for some inversion-free word $v$.

For an arbitrary word $w\in\R^\infty$, set $h_w(x):=\#\{j\dvtx w_j\leq
x\}
$. The
function $h_w\dvtx\R\to\Z_+\cup\{+\infty\}$ is weakly increasing and
right-continuous, hence it is completely determined by its restriction to
the set $\mathbb Q$ of rational
numbers.

For $x\in\R$, let $\phi_x\dvtx\R^\infty\to\{1,2\}^\infty$ be the
mapping which
replaces each $w_j\in(-\infty,x]$ by $1$ and each $w_j\in(x,+\infty)$
by $2$.
The measure $\phi_x^\infty(P)$ on $\{1,2\}^\infty$ is
$q$-exchangeable and
extreme, by virtue of Proposition \ref{extreme}. Since $h_w(x)$ is the
number of $1$'s in $\phi_x(w)$, the ergodicity implies that the value $h_w(x)$
is the same for $P$-almost all words $w$. Letting $x$ run over $\mathbb Q$,
we see that, outside a $P$-null set of words, the value $h_w(x)$ does not
depend on $w$ for each $x\in\R$; we denote by $h(x)$ this common
value. The
function $h(x)$ is again weakly increasing and right-continuous, and it assumes
values in $\Z_+\cup\{+\infty\}$.

Recall that in the $d=2$ case, $q$-exchangeability implies the
dichotomy that either
$1$ appears finitely many times and 2 appears infinitely often, or 2
does not
appear at all. From this, $h(x)\equiv\infty$ would imply $w_j\leq x$
for all
$j$, which is impossible. It follows that $h(x)$ cannot be identically
equal to
$+\infty$.

By a similar argument, $h(x)$ also cannot be identically equal to a
finite constant.

Defining
$A$ to be the set of the jump points of $h$, we see that $A$
is either a nonempty
finite set $a_1<\cdots<a_d$ or a countably infinite set of the form
$a_1<a_2<\cdots.$ In the latter case, we set $a^*=\sup\{a_i\}=\lim
a_i\in\R\cup\{+\infty\}$.
By the definition of $h(x)$, the function is constant on every interval
of the form
\[
(-\infty, a_1),\qquad [a_{i-1}, a_i),\qquad [a^*, +\infty).
\]

Finally, observe that if one ignores the $P$-null set of words mentioned above,
then any word $w$ does not contain letters from the open intervals
\[
(-\infty, a_1),\qquad (a_{i-1}, a_i),\qquad (a^*, +\infty).
\]
We conclude that $P$ is concentrated on $A^\infty$.
\begin{remark}
We note, in passing, that this argument fails for more general ordered spaces.
For instance, it cannot be applied to ${\mathbb R}^k$ $(k>1)$ with
lexicographic order because the order is not separable and $h$ cannot be
determined by its restriction to a countable set.
\end{remark}

%s8 ###
\section{Quantization}\label{8}

A motivation for studying the $q$-exchangeability is that this property
can be
viewed as a quantization of conventional exchangeability. We comment
briefly on this connection.

In the classical setting, each extreme exchangeable $P$ on $\R^\infty$
is of
the form $\nu^{\otimes\infty}$, where $\nu$ is the limit of
empirical measures,
meaning that for every Borel $B\subset\R$, as $n\to\infty$, the
random word
satisfies the strong law of large numbers
%
%e8.1 ###
%
\begin{equation}\label{LawLN}
\#\{j\leq n | w_j\in B\}\sim n \nu(B) \qquad P\mbox{-a.s.}
\end{equation}
Trivially, $0<P(w_1\in B)<1$ if and only if $0<\nu(B)<1$, in which
case letters
from $A$ appear in $w$ infinitely many times for both $A=B$ and
$A=B^c:=\R\setminus B$.

In the framework of $q$-exchangeability (with $q<1$), the analog of
(\ref{LawLN}) is
%
%e8.2 ###
%
\begin{equation}\label{LawLNq}
\#\{j\leq n | w_j\in B\}\to\nu_q(B)\qquad P\mbox{-a.s.},
\end{equation}
where $\nu_q$ is a \textit{counting} measure associated with some height function
$h$, so the letters from $B$ are represented in $w$ exactly $\nu_q(B)$
times. Similarly to the above, one sees, from the formula
\[
P(w_1\in B)=\sum_{\{x\in B | \nu_q\{x\}>0\}}
q^{\nu_q(-\infty,x)}\bigl(1-q^{\nu_q\{x\}}\bigr),
\]
that $0<P(w_1\in B)<1$ if and only if $\nu_q(B)>0$ and $\nu_q(B^c)>0$.

There are many ways to approach ex\-change\-ability via
$q$-ex\-change\-ability, that is, to obtain independent sampling in the
classical limit $q\to1$. One possible explicit realization of such a
limit is
the following quantization of homogeneous product measures.

Let $\nu$ be a probability measure on $\R$ with distribution function
$F(x):=\nu(-\infty,x]$. Let $F^{-1}(p):= \inf\{x\in\R\dvtx F(x)\geq p\}
$ be the\vspace*{2pt}
corresponding quantile function and consider the countable collection of
quantiles $\alpha_k:=F^{-1}(1-q^{k})$, $k\in\N$, as letters of the
inversion-free
word $v:=\alpha_1\alpha_2\cdots.$ The idea is to create a bridge
between independent
sampling from $\nu$ and the $q$-shuffle for the counting measure
$\nu_q=\sum_{j\in\N}\delta_{\alpha_j}$ by means of independent
sampling from
the measures
\[
\widetilde{\nu}_q=\sum_{k\in\N} G_q(k) \delta_{\alpha_k}.
\]

\begin{proposition}
As $q\to1$, for $v=\alpha_1\alpha_2\cdots,$ the $q$-shuffle
measures $P^{(v)}$
converge, in the sense of weak convergence of the finite-dimensional marginal
measures $P_n^{(v)}, n\in\N$, to the product measure $\nu^{\otimes
\infty}$.
\end{proposition}
\begin{pf}
For $\xi$ a random variable with geometric distribution $G_q$, the
distribution of randomized quantile $\alpha_\xi$ is $\widetilde{\nu
}_q$. It is
convenient to introduce two more random variables: $\zeta$ with uniform
distribution on $[0,1]$ and $\zeta_q$ with the discrete distribution
%
%e8.3 ###
%
\begin{equation}\label{meme}
\sum_{k\in\N} G_q(k) \delta_{1-q^k}.
\end{equation}
From standard properties of the quantile function, the distribution of
$F^{-1}(\zeta)$ is $\nu$ and the distribution of $F^{-1}(\zeta_q)$ is
$\widetilde{\nu}_q$, so we can identify $\alpha_\xi=F^{-1}(\zeta_q)$.

Now, the measure (\ref{meme}) was designed so that the mass of each interval
$[0,1-q^k]$ is $1-q^k$ and the largest atom has mass $1-q$, which approaches
$0$ as $q\to1$. Therefore, $\zeta_q$ converges in distribution to
$\zeta$. On
the other hand, the set of discontinuities of the quantile function is
at most
countable and so has Lebesgue measure zero, hence $F^{-1}$ preserves the
convergence relation (see, e.g., \cite{Bi}, Theorem 5.1), meaning that
$F^{-1}(\zeta_q)\to_d F^{-1}(\zeta)$. The latter is the same as
\[
{\mathbb P}(\alpha_{\xi}\leq x)\to F(x) \qquad\mbox{as }q\to1,
\]
where $x$ is an arbitrary continuity point of $F$. For any nonnegative integer
$m$, the total variation distance between $\xi$ and the shift $\xi+m$ equals
$1-q^m$, from which the above can be strengthened as
\[
{\mathbb P}(\alpha_{\xi+m}\leq x)\to F(x) \qquad\mbox{as }q\to1.
\]
Likewise, if $\xi_1,\xi_2,\ldots$ are independent copies of $\xi$ and
$m_1,\ldots,m_n$ are arbitrary fixed nonnegative integers, then we have
\[
{\mathbb P}(\alpha_{\xi_1+m_1}\leq x_1,\ldots,\alpha_{\xi
_n+m_n}\leq x_n)
\to F(x_1)\cdots F(x_n) \qquad\mbox{as }q\to1,
\]
where $x_1,\ldots,x_n$ are arbitrary continuity points of $F$.

Let $w_1w_2\cdots$ be the $q$-shuffle of $1 \cdot2 \cdots,$ constructed
from the
independent geometric $\xi_1,\xi_2,\ldots,$ as in Definition \ref{shuffle4}.
It easily follows from the definition that $\xi_j\leq w_j<\xi_j+j$,
whence the above implies
\[
{\mathbb P}(\alpha_{w_1}\leq x_1,\ldots,\alpha_{w_n}\leq x_n)
\to F(x_1)\cdots F(x_n) \qquad\mbox{as }q\to1
\]
for continuity points $x_1,\ldots, x_n$, which is precisely the property
of weak
convergence of $P^{(v)}_n$ which we wanted to prove.
\end{pf}

This construction provides quantization of homogeneous product measures
on~$\R^\infty$. Extension to the general exchangeable case is
straightforward in
the light of de Finetti's theorem: we simply randomize $\nu$.

%s9 ###
\section{Random flags over a Galois field}\label{9}

Fix $q\in(0,1)$ and set ${\tilde q}=q^{-1}$ so that ${\tilde q}>1$. In
this section, we assume
that ${\tilde q}$ is a power of a prime number.

Let $\F_{\tilde q}$ be the Galois field with ${\tilde q}$ elements and
let $V_\infty$
be an
infinite-dimensional vector space over $\F_{\tilde q}$ with a
countable basis
$\{v_1,v_2,\ldots\}$. Defining $V_n$ to be the linear span of vectors
$v_1,\ldots,v_n$, we have $\bigcup_{n\geq1} V_n=V_\infty$, so each
element of
$V_\infty$ can be uniquely written in the basis as an infinite vector with
finitely many nonzero components.

For $d\in\N$, by a \textit{decreasing $d$-flag in
$V_\infty$}, we shall mean a $(d+1)$-tuple $X=(X(i))$ of linear
subspaces in
$V_\infty$ such that
\[
V_\infty=X(0)\supseteq X(1)\supseteq\cdots\supseteq X(d-1)\supseteq
X(d)=\{0\}.
\]
Keep in mind that our definition disagrees with the conventional notion
of a
flag, in that the inclusions are not necessarily strict. In the same
way, we
define decreasing $d$-flags in each space $V_n$. Let $\XX_d(V_\infty
)$ and
$\XX_d(V_n)$ denote the sets of the decreasing $d$-flags in $V_\infty
$ and
$V_n$, respectively.
\begin{lemma}
One can identify $\XX_d(V_\infty)$ with the projective limit space
$\varprojlim\XX_d(V_n)$, where the projection $\XX_d(V_{n+1})\to\XX
_d(V_n)$ is
determined by taking the intersection with $V_n$.
\end{lemma}
\begin{pf}
Indeed, the map $\XX_d(V_\infty)\to\varprojlim\XX_d(V_n)$ is
defined by
assigning to a flag $X=(X(i))$ in $V_\infty$ the sequence
$\{X_n\in\XX_d(V_n)\}$ of flags with $X_n(i)=X(i)\cap V_n$.
Clearly, the flags
$X_n$ are consistent with the projections $\XX_d(V_{n+1})\to\XX
_d(V_n)$ and
hence determine an element of the projective limit space. The inverse map
assigns to any such sequence $\{X_n\}$ the flag $X\in\XX_d(V_\infty
)$ with
$X(i)=\bigcup X_n(i)$.
%\rightqed
\end{pf}

Using the lemma, we endow $\XX_d(V_\infty)$ with the topology of projective
limit. In other words, a small neighborhood of a flag $X=(X(i))$ is
formed by
the flags $Y=(Y(i))$ such that $X(i)\cap V_n=Y(i)\cap V_n$ for all $i$ and
some fixed large $n$. We will consider the $\sigma$-algebra of Borel
sets in $\XX_d(V_\infty)$
relative to this topology.

Let $\G_n$ be the group of all invertible linear transformations of
the space
$V_\infty$ that leave $V_n$ invariant and fix the basis vectors
$v_{n+1},v_{n+2},\ldots.$ We then have
$\{e\}=\G_0\subset\G_1\subset\G_2\subset\cdots$ and we define
$\G_\infty:=\bigcup_{n\geq1}\G_n$. The group $\G_n$ is finite and
isomorphic to
the group $\mathit{GL}(n,\F_{\tilde q})$ of invertible $n\times n$ matrices
over $\F_{\tilde q}
$. The
countable group $\G_\infty$ is isomorphic to the group $\mathit{GL}(\infty,\F
_{\tilde q}
)$ of
infinite invertible matrices $(g_{ij})$, such that $g_{ij}=\delta
_{ij}$ for
large enough $i+j$.

The group $\G_n$ acts, in a natural way, on $\XX_d(V_n)$ and the group
$\G_\infty$ acts on $\XX_d(V_\infty)$ by continuous transformations.
The next
proposition is an extension of~\cite{GO2}, Lemma 5.2.
\begin{proposition}\label{gibbs}
There exists a natural bijection $P\leftrightarrow\mathcal P$ between
$q$-exchangeable Borel probability measures on $\N_d^\infty$ and
$\G_\infty$-invariant Borel probability measures on $\XX_d(V_\infty)$.
\end{proposition}
\begin{pf}
The desired bijection is constructed by understanding $P$ as a Gibbs measure
on the path space $\Path(d)$ of the $q$-Pascal pyramid $\Ga(q,d)$,
as defined in Section \ref{5}.

We assign to $P$ a function
$\varphi(\la)$ on the vertices in the
following way. Given a vertex $\la\in\Ga(q,d)$, the probability of a finite
path ending at $\la$
equals the weight of the path times a quantity that (for given $P$)
depends only on $\la$; let us denote this quantity $\varphi(\la)$.

The Gibbs measure is uniquely determined by this function $\varphi$,
which must satisfy the rule of addition of probabilities along the path
%
%e9.1 ###
%
\begin{equation}\label{gibbs1}
\varphi(\la)=\sum_{a=1}^d \operatorname{weight}(\la,\la
+e_a)\varphi(\la+e_a)
\end{equation}
for all $\la\in\Z_+^d$, where the weight of the edge $(\la,\la
+e_w)$ is
specified in (\ref{weight}) as
%
%e9.2 ###
%
\begin{equation}\label{gibbs2}
\operatorname{weight}(\la,\la+e_a)=q^k \qquad\mbox{for } k=\la
_{a+1}+\cdots+\la_d.
\end{equation}
One must also add the normalization condition $\varphi(\bar0)=1$, which
implies that
%
%e9.3 ###
%
\begin{equation}\label{sumprob}
\sum_{\la\in\Z_+^d\dvtx |\la|=n}\dim(\la)\varphi(\la)=1,\qquad
n=1,2,\ldots,
\end{equation}
so that $\dim(\la)\varphi(\la)$ is the probability that a random walk
on $\Ga(q,d)$
driven by $P$ ever visits $\la$.

Conversely, if a nonnegative function $\varphi$ satisfies (\ref
{gibbs1}) and
the normalization condition, then it defines a Gibbs measure. Such functions
$\varphi$ play a central role in the work of Kerov and Vershik (see, e.g.,
\cite{VK}), who call them ``harmonic.'' However, this terminology is
unfortunate as it disagrees with the conventional concept of a harmonic function
in the literature on Markov processes.

We now wish to show that precisely the same functions are associated with
$\G_\infty$-invariant measures. Indeed, there is a one-to-one correspondence
between $\G_\infty$-invariant probability measures $\mathcal P$ on
$\XX_d(V_\infty)$ and sequences $\{\mathcal P_n\}$ of probability
measures such
that each $\mathcal P_n$ is a measure on $\XX_d(V_n)$, invariant under
$\G_n$,
and various $\mathcal P_n$'s are consistent with respect to the projections
$\XX_d(V_{n+1})\to\XX_d(V_n)$. Specifically, the correspondence is
established
by letting $\mathcal P_n$ be the push-forward of $\mathcal P$ under the
projection $\XX_d(V_\infty)\to\XX_d(V_n)$.

Observe that the $\G_n$-orbit of a $d$-flag $X_n=(X_n(i))\in\XX
_d(V_n)$ is
uniquely determined by the $d$-tuple of nonnegative integers
\[
\la_i=\dim V_n(i-1)-\dim V_n(i),\qquad i=1,\ldots,d,
\]
which determine a vector $\la\in\Z_+^d$ with $|\la|=n$. The reader
needs to be
warned that the dimension of a linear space over $\F_q$ in this
formula and
below should not be confused with the combinatorial dimension function
in the
Pascal pyramid, as, for instance, in (\ref{sumprob}). We will say that
the vertex $\la$
is the \textit{type} of the flag. Conversely, every such $\la$
corresponds to an
orbit. Let $\psi(\la)$ be the mass that $\mathcal P_n$ gives to each
of the
flags of type $\la$. The consistency of the measures $\mathcal P_n$ with
respect to the projections means that
%
%e9.4 ###
%
\begin{equation}\label{gibbs3}
\psi(\la)=\sum_{a=1}^d \operatorname{weight}'(\la,\la+e_a)\psi
(\la
+e_a),\qquad
\la\in\Z_+^d,
\end{equation}
where $\operatorname{weight}'(\la,\la+e_a)$ stands for the number of flags
$X_{n+1}\in\XX_d(V_{n+1})$ of type $\la+e_a$ projecting onto any
fixed flag
$X_n\in\XX_d(V_n)$ of type $\la$. Conversely, each function $\psi
(\la
)\ge0$
satisfying (\ref{gibbs3}) and the normalization condition $\psi(\bar0)=1$
determines a consistent sequence $\{\mathcal P_n\}$ and hence a
$\G_\infty$-invariant probability measure $\mathcal P$ on $\XX
_d(V_\infty)$.

We claim that
\[
\operatorname{weight}'(\la,\la+e_a)={\tilde q}^{ n-k}=q^{k-n},
\]
where $k$ is the same as in (\ref{gibbs2}), that is, $k=\dim X_n(a)$. Indeed,
if a flag $X_{n+1}$ is projected onto $X_n$, then it has type $\la+e_a$
if and
only if
\[
\dim X_{n+1}(i)=\dim X_n(i)+1 \qquad\mbox{for } 0\le i\le a-1
\]
and
\[
\dim X_{n+1}(j)=\dim
X_n(j) \qquad\mbox{for } a\le j\le d.
\]
This means that there exists a nonzero vector $v\in V_{n+1}\setminus
V_n$ such
that, for every $i=0,\ldots,a-1$, the subspace $X_{n+1}(i)$ is spanned
by $X_n(i)$
and $v$. Such a vector is uniquely defined up to a scalar multiple and addition
of an arbitrary vector from $X_n(a)$. Therefore, the number of options
is equal
to the number of lines in $V_{n+1}/X_n(a)$ not contained in
$V_n/X_n(a)$, which
equals
\[
\frac{{\tilde q}^{ n+1-k}-1}{{\tilde q}-1} - \frac{{\tilde q}^{
n-k}-1}{{\tilde q}-1}={\tilde q}^{ n-k}.
\]

Viewing equations (\ref{gibbs1}) and (\ref{gibbs3}) as recursions on
$\varphi$,
respectively, $\psi$, we see that they are similar, with the coefficients
related as
\[
\operatorname{weight}'(\la,\la+e_a)=\operatorname{weight}(\la,\la
+e_a)q^{-n},\qquad
n=|\la|.
\]
Setting
\[
\varphi(\la)=q^{n(n-1)/2}\psi(\la)
\]
yields an isomorphism $\{\varphi\}\leftrightarrow\{\psi\}$
between the convex compact sets of nonnegative solutions
to (\ref{gibbs1}) and (\ref{gibbs3}), respectively. Also, note that the
above relation does not affect the normalization condition. This
completes the proof.
\end{pf}
\begin{remark}
By virtue of the isomorphism in Proposition \ref{gibbs}, the extreme measures
$P$ correspond bijectively to extreme measures $\mathcal P$.
\end{remark}
\begin{remark}
Define a \textit{decreasing $\N$-flag} in $V_\infty$ as an infinite collection
$X=(X(i))$ of subspaces such that
\[
V_\infty=X(0)\supseteq X(1)\supseteq\cdots,\qquad \bigcap_{i\in\N}
X(i)=\{0\}.
\]
The result of Proposition \ref{gibbs} remains true when $\N_d$ is
replaced by
$\N$. That is, $q$-exchangeable probability measures on $\N^\infty$
correspond
bijectively to $\G_\infty$-invariant probability measures on the
space of
decreasing $\N$-flags. The proof is identical, except with $\Ga(q,d)$
replaced by
$\Ga(q,\infty)$.
\end{remark}
\begin{remark}
Let $V^\infty$ be the dual vector space to $V_\infty$. We endow
$V^\infty$ with
the topology of simple convergence of linear functionals; it then
becomes a
compact topological space. As an additive group, $V^\infty$ is also the
Pontryagin dual of~$V_\infty$, viewed as a discrete additive group.
Passing to the
orthogonal complement establishes a bijection between arbitrary linear
subspaces in $V_\infty$ and \textit{closed} linear subspaces in
$V^\infty$.
Define an \textit{increasing} $d$-flag in $V^\infty$ as a collection of closed
subspaces
\[
\{0\}=Y(0)\subseteq Y(1)\subseteq\cdots\subseteq Y(d)=V^\infty
\]
and an \textit{increasing} $\N$-flag in $V^\infty$ as an infinite
collection of
closed subspaces
\[
\{0\}=Y(0)\subseteq Y(1)\subseteq\cdots,\qquad \overline{
\bigcup_{i\in\N}Y(i)} =V^\infty,
\]
where the horizontal line indicates closure. By duality, the increasing
$d$-flags
in $V^\infty$ are in one-to-one correspondence with the decreasing $d$-flags
in $V_\infty$. Moreover, this correspondence is consistent with the natural
action of the group $\G_\infty$ on~$V^\infty$. The same also holds for
$\N$-flags. Thus, instead of considering invariant measures on
decreasing flags in $V_\infty$, one can equally well deal with invariant
measures on the set of increasing flags in $V^\infty$.\vadjust{\goodbreak}
\end{remark}

\begin{appendix}\label{10}
\section*{Appendix: The Mallows measure}

In this Appendix, we sketch some properties of the Mallows measures
$\Q_n$ and~$\Q$.
To state the results, we need some preparation. It is convenient to
represent a
generic permutation $\si\in\Sym_n$ as an $n\times n$ permutation matrix
$\si(i,j)$, where the entry $\si(i,j)$ equals 1 or 0, depending on
whether or not
$\si(j)=i$.
Such permutation matrices are \textit{strictly monomial}, in the sense
that they have
one and only one nonzero element per row and per column.
Note
that this realization of permutations by strictly monomial matrices
takes the
group multiplication into conventional matrix multiplication and the
inversion map $\si\mapsto\si^{-1}$ corresponds to matrix transposition.
Likewise, the group $\Sym$ can be realized as the group of strictly monomial
matrices of infinite size.

More generally, a 0--1 matrix of finite or infinite size is \textit{weakly
monomial} if each row and each column contains \textit{at most} one 1, the
other entries being $0$'s. Let $M(n)$ and $M$ denote the sets of weakly
monomial 0--1 matrices of size $n\times n$ and $\infty\times\infty$,
respectively.
Both $M(n)$ and $M$ are semigroups under matrix multiplication and
$\Sym_n\subset M(n)$ and $\Sym\subset M$ are respective subgroups of
invertible
elements. An additional operation in $M(n)$ and $M$ is matrix
transposition, which is an involutive antiautomorphism.

For $k=1,2,\ldots,$ the \textit{truncation operation} $\tth_k$ assigns
to a
matrix of size $\infty\times\infty$ or $l\times l$ with $l\geq k$ the
$k\times
k$ submatrix comprised of the entries $(i,j)$ with $i,j\le k$. Obviously,
$\tth_k$ projects $M(n)$ onto $M(k)$ for any $n>k$. Likewise, $\tth
_k$ projects
$M$ onto $M(k)$. Using these projections, we may identify $M$ with the
projective limit space $\varprojlim M(k)$. We endow $M$ with the corresponding
projective limit topology; $M$ then becomes a compact topological
space. By definition, a fundamental system of neighborhoods of a matrix
$m\in M$ is
formed by the subsets $\{m'\in M\mid\tth_k(m')=\tth_k(m)\}$,
$k=1,2,\ldots.$

It is readily checked that the restriction of $\tth_k\dvtx M\to M(k)$ to
the subset
$\Sym\subset M$ is surjective for every $k$. It follows that $\Sym$ is
dense in
$M$ (and even $\Sym_\infty$ is dense). Recall that we have endowed
$\Sym
$ with the
$\sigma$-algebra of Borel sets inherited via the embedding $\Sym
\subset
\N^\infty$.
Clearly, this Borel structure coincides with that induced by the embedding
$\Sym\subset M$. Thus, any Borel probability measure on $\Sym$ or on
$\Sym_n\subset\Sym$ can be viewed as a measure on $M$ (here, we identify
$\Sym_n$ with the subgroup in $\Sym$ fixing all integers from
$\N\setminus\N_n$). In particular, we may view the Mallows measures
$\Q
_n$ and
$\Q$ as probability measures on the compact space $M$. This makes sense
of the
following assertion.
\setcounter{definition}{0}
\begin{proposition}\label{mallows1}
As $n\to\infty$, $\Q_n$ weakly converge to $\Q$.
\end{proposition}
\begin{pf}
Let $\tth_k(\Q_n)$ and $\tth_k(\Q)$ denote the respective push-forwards
of $\Q_n$ and $\Q$
under $\tth_k$. By the definition of the topology on $M$ and the finiteness
of $M(k)$, it suffices to prove that for any $k$ and any fixed matrix
$m\in
M(k)$, $\tth_k(\Q_n)(\{m\})$ converges to $\tth_k(\Q)(\{m\})$.

Taking into account Remark \ref{inversion}, it is convenient to replace
$\Q_n$
and $\Q$ by their respective push-forwards under the matrix
transposition; let us denote
them as $\Q_n'$ and $\Q'$, respectively. Thus, we will prove the equivalent
assertion that $\tth_k(\Q'_n)(\{m\})$ converge to $\tth_k(\Q)(\{m\})$.

Let $w=w_1w_2\cdots$ be the output of the $q$-shuffling algorithm
applied to the
infinite word $1\cdot2\cdots.$ As usual, we identify $w$ with the random permutation
$\si\in\Sym$ by writing $w=\si(1)\si(2)\cdots.$ From this,
one sees
that the
quantity $\tth_k(\Q')(\{m\})$ is equal to the probability of the event
that for
each $j=1,\ldots, k$, the letter $w_j$ either equals some $i\in\{
1,\ldots,k\}$ if the
matrix $m$ has 1 in the $j$th column in position $(i,j)$,
or $w_j>k$ if the $j$th column of $m$ consists entirely of $0$'s.

For instance, if $m=\bigl[{ 0 \atop0} \enskip
{1 \atop0}\bigr]\in M(2)$, then the event in question is that the
first step
of the algorithm yields $w_1>2$ and the second step yields $w_2=1$.

The quantity $\tth_k(\Q'_n)(\{m\})$ admits exactly the same
interpretation in
terms of the finite $q$-shuffle applied to the finite word $1\cdots n$.

Now, the desired convergence of the probabilities follows from the fact that
as $n\to\infty$, the truncated geometric distributions directing the finite
$q$-shuffle (Definition \ref{shuffle2}) converge to the infinite geometric
distribution directing the infinite $q$-shuffle (Definition \ref{shuffle4}).
\end{pf}
\begin{corollary}\label{symmetry}
The Mallows measures $\Q_n$ and $\Q$ are invariant under the
group inversion map $\si\mapsto\si^{-1}$.
\end{corollary}
\begin{pf}
Given a matrix $m\in M(n)$, let us say that two distinct positions
$\{(i_1,j_1), (i_2,j_2)\}$ occupied by 1's are \textit{in inversion} if
the two
differences $i_1-i_2$ and $j_1-j_2$ have opposite signs (note that these
differences cannot vanish) and denote by $\inv(m)$ the total number of
unordered pairs of positions in inversion. Clearly, $\inv(m)=\inv
(m')$, where
$m'$ stands for the transposed matrix.

On the other hand, if $\si\in\Sym_n$ and $m:=[\si(i,j)]$ is the
corresponding
permutation matrix, then we obviously have $\inv(\si)=\inv(m)$. If
$\si
$ is
replaced by $\si^{-1}$, then $m$ is replaced by $m'$. Therefore,
$\inv(\si)=\inv(\si^{-1})$, which implies the desired symmetry property
of $\Q_n$.
The analogous property for $\Q$ now follows from Proposition \ref{mallows1}.
\end{pf}
\begin{remark}\label{symmetry0}
The ``absorption sampling'' mentioned above (see \cite{Ke2} for
history and
references) seems not to have been identified with the Mallows measure
on $M$.
This connection, along with the invariance of $\Q$ under matrix
transposition, make obvious the unexplained symmetry in formulae like
\cite{Ke1}, equation (10) and~\cite{B}, equation (2.12).

Likewise, the number of inversions is also invariant under reflection
with respect to the secondary matrix diagonal, which swaps $(i,j)$ and
$(n+1-j,n+1-i)$, so $\Q_n$ is also preserved by this transformation.
However, this operation has no analog for the infinite group $\Sym$.
\end{remark}
\begin{remark}\label{symmetry2}
Observe that the group $\Sym_\infty$ acts on $\Sym$ both by left and right
shifts: an element $\si\in\Sym_\infty$ maps an element $\tau\in
\Sym$ to
$\si\tau$ or $\tau\si^{-1}$, respectively. Under the right action, the
elementary transposition $\si_i:=(i,i+1)\in\Sym_\infty$ swaps the
letters of a
permutation word $\wt\tau$ in the $i$th and $(i+1)$th positions,
while under
the left action, the same element $\si_i$ swaps the letters $i$ and
$(i+1)$ in
$\wt\tau$. That is, under the right action on permutation words, we look
at positions, while under the left action, we look at the letters themselves.
The inversion map intertwines both actions.

We know that $\Q$ is a unique probability measure on $\Sym$ that is
quasi-invariant under the right action, with a special cocycle, (\ref{rho-q}).
The symmetry property of the measure $\Q$ implies that it is also
quasi-invariant
under the left action. To compute the corresponding cocycle, we return to
the definition (\ref{stable}) of the additive cocycle and observe that instead
of taking the $n$-truncated word with large $n$, we can equally well
deal with
arbitrary finite subwords, provided that they are large enough. Using this
reformulation, we see that the additive cocycle is preserved under the group
inversion on $\Sym$, as is the corresponding multiplicative cocycle.

It follows that the cocycle corresponding to the left action remains
the same.
Consequently, $\Q$ can also be characterized as a unique probability
measure on
$\Sym$ which is quasi-invariant under the left action of $\Sym_\infty
$ with
the same cocycle as before.
\end{remark}

The next proposition describes the finite-dimensional distributions of the
Mallows measure $\Q$ viewed as a measure on $M=\varprojlim M(k)$. We
use the
following notation: $m$ is an arbitrary matrix from $M(k)$;
$I\subset\{1,\ldots,k\}$ is the set of indices of the rows in $m$
containing 1's;
$J\subset\{1,\ldots,k\}$ is the set of indices of the columns in $m$ containing
1's; $r=|I|=|J|$ is the rank of $m$; $\inv(m)$ has the same meaning as
in the
proof of Corollary \ref{symmetry}.
\begin{proposition}\label{mallows3}
Using the above notation,
%
%e9.5 ###
%
\setcounter{equation}{0}
\begin{equation}\label{mallows4}
\tth_k(\Q)(\{m\})=(1-q)^{r}q^{k^2-2kr-r+\inv(m)+\sum_{i\in I}i+\sum
_{j\in J}j}.
\end{equation}
\end{proposition}
\begin{pf}
We apply the same method as in Section \ref{6}, that is, reduce the
alphabet $\N$ to the
finite alphabet $\N_{k+1}$ using the monotone map $f_{k+1}(a)=a\wedge(k+1)$.
The key idea is that if $w=w_1w_2\cdots=\si(1)\si(2)\cdots$ is the
random output
of the infinite $q$-shuffle of the word $v=1 \cdot2\cdots$ then, as seen
from the
proof of Proposition \ref{mallows1}, the truncated matrix $\tth_k(\si
)$ depends
only on the first $k$ letters of the word $f_{k+1}^\infty(w)$ (i.e.,
all of the letters $\ge k+1$ become indistinguishable).

On the other hand, by virtue of Proposition \ref{monotone}, the random word
$f_{k+1}^\infty(w)$ is the output of the infinite $q$-shuffle applied
to the
inversion-free word
\[
v':=1\cdots k \underbrace{(k+1)(k+1)\cdots}_{\infty}\in(\N
_{k+1})^\infty.
\]
In the notation of Section \ref{4}, the law of the random word
$f_{k+1}^\infty(w)$ is given by the measure $P^{(v')}$ and the
distribution of
the first $k$ letters is given by the marginal $P^{(v')}_k$, for which we
have an explicit expression; see (\ref{marg1}). In this formula, we
need to take
\begin{eqnarray*}
l_1&=&\cdots=l_k=1,\qquad l_{k+1}=\infty,\qquad \mu_{k+1}=k-r,
\\
\mu_a&=&\cases{1, &\quad $a\in I$,\cr0, &\quad $a\in\{1,\ldots,k\}\setminus I$,}
\end{eqnarray*}
and then the direct computation gives (\ref{mallows4}).
\end{pf}

There is another way of approximating $\Q$ by the $\Q_n$'s. Namely, we
will see
that $\Q$ can be represented as the projective limit of the $\Q_n$'s.
Incidentally, we will realize $\Q$ as a product measure.

As usual, we will identify permutations with the corresponding permutation
words. For any $n\ge2$, we define the projection $\Sym_n\to\Sym
_{n-1}$ as
the deletion of $n$ from a permutation word. Using these projections,
we construct
the projective limit space $\varprojlim\Sym_n$, which is a compact topological
space in the standard topology. We have a natural embedding
%
%e9.6 ###
%
\begin{equation}\label{embed}
\Sym\hookrightarrow\varprojlim\Sym_n ,
\end{equation}
which is specified by the projection $\Sym\to\Sym_n$ which removes
all
letters larger than $n$ from an infinite permutation word.

Note that $\Sym$ is a \textit{proper} subset of $\varprojlim\Sym_n$. Indeed,
there is a natural one-to-one correspondence between elements of
$\varprojlim\Sym_n$ and all possible linear orders on the set $\N$,
of which
the orders induced by permutation words $\si(1)\si(2)\cdots$
comprise a
relatively small part. Still, $\Sym$ is dense in $\varprojlim\Sym_n$.
\begin{proposition}
The measures $\Q_n$ are consistent with the projections $\Sym_n\to
\Sym_{n-1}$,
so we can define the projective limit $\Q_\infty:=\varprojlim\Q_n$,
which is a probability measure on $\varprojlim\Sym_n$. The image of
$\Sym$
under the embedding (\ref{embed}) has full $\Q_\infty$-measure and the
restriction of $\Q_\infty$ to $\Sym$ coincides with the Mallows
measure~$\Q$.
\end{proposition}
\begin{pf} For a permutation $\si\in\Sym_n$ (which we identify with the
corresponding permutation word), set
\[
\wt\be_j=\wt\be_j(\si)=\#\{i<j\mid\mbox{$i$ precedes $j$}\}+1,
\qquad j=1,\ldots,n
\]
[cf. (\ref{backrank})]. The link with (\ref{backrank}) is the identity
$\wt\be_j(\si)=\be_j(\si^{-1})$.\vspace*{2pt}

The correspondence $\si\mapsto(\wt\be_1(\si),\ldots,\wt\be_n(\si
))$ is a
bijection,
%
%e9.7 ###
%
\begin{equation}\label{bijection}
\Sym_n\to\N_1\times\cdots\times\N_n,
\end{equation}
and we have a counterpart of Proposition \ref{ranks-fin}: under $\Q
_n$, the
coordinates $\wt\be_j$ are independent and $j+1-\wt\be_j$ is distributed
according to $G_{q,j}$. This can be deduced from Proposition \ref{ranks-fin}
taken together with the symmetry property of $\Q_n$ (Proposition
\ref{symmetry}), or can be easily checked directly.

Under the bijection (\ref{bijection}), the projection $\Sym_n\to\Sym
_{n-1}$ is
simply the deletion of the last letter. This enables us to identify
$\varprojlim\Sym_n$ with the infinite product space $\prod
_{n=1}^\infty
\N_n$. Under
this identification, the measure $\varprojlim\Q_n$ becomes the
product of
truncated geometric distributions. The image of $\Sym$ in
$\prod_{n=1}^\infty\N_n$ consists of those sequences $(i_1,i_2,\ldots
)$ for
which $i_n\to\infty$. {}From this, it is readily checked that $\Sym$
has full
measure.

It remains to check that the measure $\varprojlim\Q_n$ coincides on
$\Sym$ with the measure~$\Q$. To this end, we use the characterization
of $\Q$ in terms of the left action of $\Sym_\infty$, as described in
Remark \ref{symmetry2}. It is easy to see that the measure
$\varprojlim\Q_n$ has the same transformation property with respect to
the left action of elementary transpositions~$\si_i$. Consequently,
$\varprojlim\Q_n=\Q$.
\end{pf}

Alternatively, one can use another chain of projections, such that the
projection $\Sym_n\to\Sym_{n-1}$ first cuts the last letter in
$\sigma(1)\cdots\sigma(n)$, then relabels the letters
$\sigma(1)\cdots\sigma(n-1)$ by the increasing bijection with $\N
_{n-1}$. A
random element of $\Sym$ under $\Q$ is representable by an infinite
sequence of
backward ranks $(\beta(1),\beta(2),\ldots)$, which are independent
and have
distribution as in Proposition~\ref{ranks-fin}.
\end{appendix}

\section*{Acknowledgments}
We are indebted to Yuliy Baryshnikov and Persi Diaconis for
illuminating discussions and references.

% imsref loaded by lrinkeviciute, 2010-06-23 10:45:16
%

%
\printaddresses


\begin{thebibliography}{19}

%b1 ###
\bibitem{Al}
%
\begin{bincollection}[mr]
\bauthor{\bsnm{Aldous},~\bfnm{David~J.}\binits{D.~J.}}
(\byear{1985}).
\btitle{Exchangeability and related topics}.
In \bbooktitle{\'{E}cole D'\'et\'e de Probabilit\'es de {S}aint-{F}lour,
{XIII}---1983}.
\bseries{Lecture Notes in Math.}
\bvolume{1117}
\bpages{1--198}.
\bpublisher{Springer}, \baddress{Berlin}.
\bid{mr={883646}}
\end{bincollection}
%
\endbibitem

%b2 ###
\bibitem{An}
%
\begin{bbook}[vtex]
\bauthor{\bsnm{Andrews},~\bfnm{George~E.}\binits{G.~E.}}
(\byear{1998}).
\btitle{The Theory of Partitions}.
\bpublisher{Cambridge Univ. Press}, \baddress{Cambridge}.
\bid{mr={1634067}}
\end{bbook}
%
\endbibitem

%b3 ###
\bibitem{B}
%
\begin{barticle}[mr]
\bauthor{\bsnm{Barakat},~\bfnm{Richard}\binits{R.}}
(\byear{1985}).
\btitle{Probabilistic aspects of particles transiting a trapping
field: An
exact combinatorial solution in terms of {G}auss polynomials}.
\bjournal{Z. Angew. Math. Phys.}
\bvolume{36}
\bpages{422--432}.
\bid{doi={10.1007/BF00944633}, mr={797237}}
\end{barticle}
%
\endbibitem

%b4 ###
\bibitem{BD}
%
\begin{barticle}[mr]
\bauthor{\bsnm{Bayer},~\bfnm{Dave}\binits{D.}} \AND
\bauthor{\bsnm{Diaconis},~\bfnm{Persi}\binits{P.}}
(\byear{1992}).
\btitle{Trailing the dovetail shuffle to its lair}.
\bjournal{Ann. Appl. Probab.}
\bvolume{2}
\bpages{294--313}.
\bid{mr={1161056}}
\end{barticle}
%
\endbibitem

%b5 ###
\bibitem{Benjamini}
%
\begin{barticle}[mr]
\bauthor{\bsnm{Benjamini},~\bfnm{Itai}\binits{I.}},
\bauthor{\bsnm{Berger},~\bfnm{Noam}\binits{N.}},
\bauthor{\bsnm{Hoffman},~\bfnm{Christopher}\binits{C.}} \AND
\bauthor{\bsnm{Mossel},~\bfnm{Elchanan}\binits{E.}}
(\byear{2005}).
\btitle{Mixing times of the biased card shuffling and the asymmetric exclusion
process}.
\bjournal{Trans. Amer. Math. Soc.}
\bvolume{357}
\bpages{3013--3029 (electronic)}.
\bid{doi={10.1090/S0002-9947-05-03610-X}, mr={2135733}}
\end{barticle}
%
\endbibitem

%b6 ###
\bibitem{Bi}
%
\begin{bbook}[vtex]
\bauthor{\bsnm{Billingsley},~\bfnm{Patrick}\binits{P.}}
(\byear{1999}).
\btitle{Convergence of Probability Measures},
\bedition{2nd} ed.
%Statistics}.
\bpublisher{Wiley}, \baddress{New York}.
\bid{doi={10.1002/9780470316962}, mr={1700749}}
\end{bbook}
%
\endbibitem

%b7 ###
\bibitem{DiaconisRam}
%
\begin{barticle}[mr]
\bauthor{\bsnm{Diaconis},~\bfnm{Persi}\binits{P.}} \AND
\bauthor{\bsnm{Ram},~\bfnm{Arun}\binits{A.}}
(\byear{2000}).
\btitle{Analysis of systematic scan {M}etropolis algorithms using
{I}wahori--{H}ecke algebra techniques}.
\bjournal{Michigan Math. J.}
\bvolume{48}
\bpages{157--190}.
\bid{doi={10.1307/mmj/1030132713}, mr={1786485}}
\end{barticle}
%
\endbibitem

%b8 ###
\bibitem{GO1}
%
\begin{barticle}[mr]
\bauthor{\bsnm{Gnedin},~\bfnm{Alexander}\binits{A.}} \AND
\bauthor{\bsnm{Olshanski},~\bfnm{Grigori}\binits{G.}}
(\byear{2006}).
\btitle{The boundary of the {E}ulerian number triangle}.
\bjournal{Mosc. Math. J.}
\bvolume{6}
\bpages{461--475, 587}.
\bid{mr={2274860}}
\end{barticle}
%
\endbibitem

%b9 ###
\bibitem{GO2}
%
\begin{barticle}[vtex]
\bauthor{\bsnm{Gnedin},~\bfnm{Alexander}\binits{A.}} \AND
\bauthor{\bsnm{Olshanski},~\bfnm{Grigori}\binits{G.}}
(\byear{2009}).
\btitle{A {$q$}-analogue of de {F}inetti's theorem}.
\bjournal{Electron. J. Combin.}
\bvolume{16}
\bpages{R78}.
\bid{mr={2529787}}
\end{barticle}
%
\endbibitem

%b10 ###
\bibitem{GS}
%
\begin{barticle}[vtex]
\bauthor{\bsnm{Greschonig},~\bfnm{Gernot}\binits{G.}} \AND
\bauthor{\bsnm{Schmidt},~\bfnm{Klaus}\binits{K.}}
(\byear{2000}).
\btitle{Ergodic decomposition of quasi-invariant probability measures}.
\bjournal{Colloq. Math.}
\bvolume{84/85}
\bpages{495--514}.
\bid{mr={1784210}}
\end{barticle}
%
\endbibitem

%b11 ###
\bibitem{K}
%
\begin{bbook}[vtex]
\bauthor{\bsnm{Kallenberg},~\bfnm{Olav}\binits{O.}}
(\byear{2005}).
\btitle{Probabilistic Symmetries and Invariance Principles}.
\bpublisher{Springer}, \baddress{New York}.
\bid{mr={2161313}}
\end{bbook}
%
\endbibitem

%b12 ###
\bibitem{Ke1}
%
\begin{barticle}[mr]
\bauthor{\bsnm{Kemp},~\bfnm{Adrienne~W.}\binits{A.~W.}}
(\byear{1998}).
\btitle{Absorption sampling and the absorption distribution}.
\bjournal{J. Appl. Probab.}
\bvolume{35}
\bpages{489--494}.
\bid{mr={1641849}}
\end{barticle}
%
\endbibitem

%b13 ###
\bibitem{Ke2}
%
\begin{bincollection}[vtex]
\bauthor{\bsnm{Kemp},~\bfnm{A.~W.}\binits{A.~W.}}
(\byear{2001}).
\btitle{A characterization of a distribution
arising from absorption sampling}.
In \bbooktitle{Probability and Statistical Models and Applications}
(\beditor{C. A. Charalambides et al.}, eds.)
\bpages{239--246}.
\bpublisher{Chapman and Hall/CRC Press}, \baddress{Boca Raton, FL}.
\end{bincollection}
%
\endbibitem

%b14 ###
\bibitem{KOO}
%
\begin{barticle}[vtex]
\bauthor{\bsnm{Kerov},~\bfnm{Sergei}\binits{S.}},
\bauthor{\bsnm{Okounkov},~\bfnm{Andrei}\binits{A.}} \AND
\bauthor{\bsnm{Olshanski},~\bfnm{Grigori}\binits{G.}}
(\byear{1998}).
\btitle{The boundary of the {Y}oung graph with {J}ack edge multiplicities}.
\bjournal{Int. Math. Res. Not.}
\bvolume{1998}
\bpages{173--199}.
\bid{doi={10.1155/S1073792898000154}, mr={1609628}}
\end{barticle}
%
\endbibitem

%b15 ###
\bibitem{Mallows}
%
\begin{barticle}[mr]
\bauthor{\bsnm{Mallows},~\bfnm{C.~L.}\binits{C.~L.}}
(\byear{1957}).
\btitle{Non-null ranking models. {I}}.
\bjournal{Biometrika}
\bvolume{44}
\bpages{114--130}.
\bid{mr={0087267}}
\end{barticle}
%
\endbibitem

%b16 ###
\bibitem{Rawlings}
%
\begin{barticle}[mr]
\bauthor{\bsnm{Rawlings},~\bfnm{Don}\binits{D.}}
(\byear{1997}).
\btitle{Absorption processes: Models for {$q$}-identities}.
\bjournal{Adv. in Appl. Math.}
\bvolume{18}
\bpages{133--148}.
\bid{doi={10.1006/aama.1996.0504}, mr={1430385}}
\end{barticle}
%
\endbibitem

%b17 ###
\bibitem{St}
%
\begin{barticle}[vtex]
\bauthor{\bsnm{Stanley},~\bfnm{Richard~P.}\binits{R.~P.}}
(\byear{2001}).
\btitle{Generalized riffle shuffles and quasisymmetric functions}.
\bjournal{Ann. Comb.}
\bvolume{5}
\bpages{479--491}.
\bid{doi={10.1007/s00026-001-8023-7}, mr={1897637}}
\end{barticle}
%
\endbibitem

%b18 ###
\bibitem{Starr}
%
\begin{barticle}[mr]
\bauthor{\bsnm{Starr},~\bfnm{Shannon}\binits{S.}}
(\byear{2009}).
\btitle{Thermodynamic limit for the {M}allows model on {$S\sb n$}}.
\bjournal{J. Math. Phys.}
\bvolume{50}
\bpages{095208, 15}.
\bid{doi={10.1063/1.3156746}, mr={2566888}}
\end{barticle}
%
\endbibitem

%b19 ###
\bibitem{VK}
%
\begin{barticle}[vtex]
\bauthor{\bsnm{Vershik},~\bfnm{A.~M.}\binits{A.~M.}} \AND
\bauthor{\bsnm{Kerov},~\bfnm{S.~V.}\binits{S.~V.}}
(\byear{1987}).
\btitle{Locally semisimple algebras: Combinatorial theory and the
$K_0$-functor}.
\bjournal{J. Math. Sci. (N. Y.)}
\bvolume{38}
\bpages{1701--1733}.
\end{barticle}
%
\endbibitem

\end{thebibliography}
\end{document}